\documentclass{article}
\usepackage{latexsym,amssymb}
\usepackage{psfrag}
\usepackage[dvips]{graphicx}

\newcommand{\C}{\mathbb{C}}

\newcommand{\R}{\mathbb{R}}
\newcommand{\T}{\mathbb{T}}
\newcommand{\<}{\langle}
\newcommand{\m}{\rangle}
\newcommand{\PP}{\mathcal{P}}
\newcommand{\G}{\mathcal{G}}
\newcommand{\f}{\mathcal{F}}

\newcommand{\p}{\mathcal{P}}
\newcommand{\Fix}{\mathop{\rm Fix}}
\newcommand{\Ker}{\mathop{\rm Ker}}
\newcommand{\CoAd}{\mathop{\rm CoAd}\nolimits}
\newcommand{\Ad}{\mathop{\rm Ad}\nolimits}
\newcommand{\so}{\mathfrak{so}(3)}
\newcommand{\ep}{\Box}

\def \Z{{\sf Z\hspace{-.45em}Z}}
\def\proof{\noindent\textbf{Proof. }}
\def\rmk{\noindent\textbf{Remark. }}
\def\rmks{\noindent\textbf{Remarks. }}
\def\ep{$\Box$\\}

\def\tet{\theta}

\def\rum{relative equilibrium }
\def\rumb{relative equilibrium}
\def\ria{relative equilibria }
\def\riab{relative equilibria}
\newtheorem{theo}{Theorem}[section]
\newtheorem{prop}[theo]{Proposition}
\newtheorem{define}[theo]{Definition}

\begin{document}

\noindent{\Large\textbf{Point vortices on the sphere:
a case with opposite vorticities}}\\
\\
F Laurent-Polz\\
Institut Non Lin\'eaire de Nice, Universit\'e de Nice Sophia Antipolis, Valbonne, France\\

\noindent\textbf{Abstract.}
We study  systems formed of $2N$ point vortices on a sphere with $N$ vortices of strength $+1$ and $N$ vortices of strength $-1$.
In this case, the  Hamiltonian is conserved by the symmetry which exchanges the positive vortices with the negative vortices.
We prove the existence of some fixed and relative equilibria, and then study their stability with the ``Energy Momentum Method''.
Most of the results obtained are nonlinear stability results.     
To end, some bifurcations are described.\\


\noindent AMS classification scheme number: 58Z05, 70H14, 70H33\\

\section{Introduction}
        
The study of the motion of many eddies (\emph{vortices}) started one century ago when Helmhotz \cite{H} introduced the point vortex model. 
With progress in dynamical systems in recent decades, the $N$-vortex problem has re-appeared.
Numerous papers have been written on vortices in the plane \cite{Ar82,Ar83,Ar83',Ar98,LR96,Sa92}.
From a geostrophic point of view, it is interesting to study point vortices on the sphere.
In this paper we consider a non-rotating sphere, that is we focus only on the spherical geometry of the Earth. 
In a forthcoming paper we  will consider a rotating sphere and then deal with the \emph{forced symmetry breaking} theory \cite{CL00}.
As many results in this paper depend only on the phase space and the symmetries of the Hamiltonian, the study of point vortices on the sphere may also be fruitful for other problems which arise in physics: 
charged particles on a sphere, $\alpha$-Euler model in turbulence, motion of cells on a spherical flame, zeros of wavefunctions in quantum physics\dots

The motion of three vortices of arbitrary strength on a sphere is studied in \cite{KN98} and \cite{KN00}. 
The stability of the relative equilibria given in \cite{KN98} is computed in \cite{PM98} and numerical simulations are done in \cite{MPS99}.
Relative equilibria formed of $N$ vortices are treated in \cite{LMR00} and the stability of a latitudinal ring of $N$ identical vorticities is computed in \cite{LMR}:
it appears that linear stability results obtained by \cite{PD93} are in fact Lyapunov stability results. 

The dynamics of point vortices on a sphere is given by an Hamiltonian system (\cite{B77}).
The Hamiltonian is invariant under rotations of the sphere, reflections of the sphere and permutations of identical vortices.
In this paper, we consider $N$ vortices of vorticity $+1$ and $N$ vortices of vorticity $-1$, where one has a new symmetry for the Hamiltonian:
the permutation which exchanges the positive vortices with the negative vortices.
However, some of these symmetries (eg reflections) are not symmetries of the equations of motion, they are \emph{time-reversing} (or, \emph{anti-symplectic}) symmetries  (Section \ref{symsection}).
From Noether's theorem, the rotational symmetry provides three conserved quantities, the components of the \emph{centre of vorticity} map $\Phi: \p \to \R^3$ where $\p$ is the phase space.

In Section 3, we determine some equilibria and \emph{\riab}.
Relative equilibria are orbits of the group action which are invariant under the flow and correspond of the motions of the point vortices which are stationary in a steadily rotating frame.
In the same way in which equilibria are critical points of the Hamiltonian $H$, \ria are critical points of the restrictions of $H$ to the level sets $\Phi^{-1}(\mu)$.
The fixed point set of a group of symmetries with some anti-symplectic symmetries might not be invariant under the flow contrary to those formed of symplectic symmetries.
However, it follows from the \emph{Principle of Symmetric Criticality} \cite{P79,Mi71} that critical points of the restriction of $H$ to such a fixed point set are critical points of the full Hamiltonian, and are therefore equilibrium points.
An analogous statement holds for \riab.
Relative equilibria derived with these results do not depend on the expression for the Hamiltonian, they depend only on the symmetries, the phase space and the momentum map.
This technique is especially useful for finding \ria in high dimensional systems with large symmetry groups.
It is probable that there are ``asymmetric'' \ria but they can not be found by the symmetry adapted methods of this paper. 
See \cite{Ar98} for examples of asymmetric relative equilibria in the plane.

In \cite{LMR00}, the authors consider also anti-symplectic symmetries and use the Principle of Symmetric Criticality, but the approach is quite different.
In the case of $N$ identical vortices, the Hamiltonian goes to minus infinity when some vortices collide.
Thus the Hamiltonian must have a maximum somewhere, so there exists an equilibrium point.
They generalized this idea to fixed point sets (see the remark after Theorem \ref{interm} for an example).

Section 4 is devoted to analyzing the stability of the \ria determined in Section 3.
The suitable concept of stability for \ria in Hamiltonian system is Lyapunov stability \emph{modulo a subgroup}.
The stability study is realized using the \emph{energy momentum method} \cite{Si91,Ma92,Pa92,LS98,OR99}.  
This is stronger than linear stability.  
For configurations with a small number of vortices, there are some Lyapunov stable relative equilibria.
For a large number of vortices  the relative equilibria we study are all linearly unstable, and in general we do not know whether there are other stable relative equilibria.
When a \rum is claimed to be linearly stable, this means in addition that the Lyapunov stability criterion (the \emph{energy momentum theorem}, a kind of Dirichlet's criterion for \riab) does not apply.
Finally, we describe some pitchfork and Hamiltonian-Hopf bifurcations in systems of four and six vortices.

\section{The point vortices system on the sphere}

\subsection{Equations of motion}

Let $N$ be an integer, $N\geq 2$.
The phase space for the $2N$-vortex problem is a product of $N$ spheres where the diagonal has been removed in order to avoid collisions:
$$
\p = \lbrace (x_1,\dots,x_{2N}) \in S^2 \times \cdots \times S^2 \mid  x_i \neq x_j\ \mbox{if}\ i \neq j \rbrace.$$
In the previous expression, each $x_i$ represents a point vortex on the sphere, of vorticity $\lambda_i$.

In principle we assume that half of the vortices are of vorticity $+1$ and the other half of vorticity $-1$, though this assumption can be relaxed when polar vortices are considered.
We could write an element of $\cal{P}$ in two different manners:

\noindent - either $x=(x_1,\dots,x_{N},x_{N+1},\dots,x_{2N})$ with $\lambda_j=+1$ and $\lambda_{N+j}=-1$, $j=1,\dots,N$, 

\noindent - or $x=(x_1,\dots,x_{N},x_1^{\prime},\dots,x_{N}^{\prime})$  with $\lambda_j=+1$ and $\lambda_j^{\prime}=-1$, $j=1,\dots,N$.
   
For notational convenience, arrows above vectors will be omitted most of the time.
We call a \textit{($+$)vortice} (resp. a \textit{($-$)vortice}) a vortice with vorticity $+1$ (resp. $-1$), a \textit{($+$)ring} (resp. a \textit{($-$)ring}) a latitudinal regular gon of ($+$)vortices (resp. ($-$)vortices), and a \textit{($\pm$)ring} a latitudinal regular polygon formed of alternating ($+$)vortices and ($-$)vortices.
Throughout this paper, $r_\alpha$ denotes the rotation of $SO(3)$ about the $(Oz)$ axis with angle $\alpha$, and $s_x$, $s_y$, $s_z$ denote respectively the reflections of $O(3)$ of planes $(yOz)$, $(xOz)$, $(xOy)$.

The equations governing the motion of $2N$ vortices are given by the following system of  $6N$ equations \cite{B77,KN98}:
$$
\dot {x_i} = X_H(x)_i=\sum_{j,j\neq i} \frac{\lambda_{j}(x_{j}\times x_{i})}{1-x_{i}\cdot x_{j}}\  ,\ i=1,\dots,2N .
$$
Since the vortices are constrained to lie on the sphere, these equations are not independent.

Let $\theta_i$ be the  co-latitude and $\phi_i$ be the  longitude of $x_i$. 
In spherical coordinates, the  system is formed of the following $4N$ equations \cite{KN98}:
$$
\dot\theta_i=-2\sum_{j=1,j\neq i}^{2N} \lambda_j\frac{\sin\theta_j\sin(\phi_i-\phi_j)}{l_{ij}^2}\ \ ,\ i=1,\dots,2N
$$
$$
\sin\theta_i\ \dot \phi_i=2\sum_{j=1,j\neq i}^{2N} \lambda_j \frac{\sin\theta_i \cos\theta_j -\sin\theta_j \cos\theta_i \cos(\phi_i-\phi_j)} {l_{ij}^2}\ \ ,\ i=1,\dots,2N
$$
where $l_{ij}^2=2(1-\cos\theta_i \cos\theta_j -\sin\theta_i \sin\theta_j \cos(\phi_i-\phi_j))$ is the square of the Euclidian distance $\|x_i-x_j\|$.
Introducing the conjugate variables $q_i=\sqrt{\vert\lambda_i\vert}\cos\theta_i$ and $p_i=sign(\lambda_i)\sqrt{\vert\lambda_i\vert}\phi_i$, and the Hamiltonian 
\begin{equation}
\label{HAM}
H = \sum_{i<j} \lambda_i  \lambda_j \ln l_{ij}^2
\end{equation}
we put the system in Hamiltonian canonical form:
$$
\dot q_i=\frac{\partial H}{\partial p_i}
$$
$$
\dot p_i=-\frac{\partial H}{\partial q_i}
$$
for all $i=1,\dots,2N$.

The phase space has so a symplectic  structure given by
$$\omega=\sum_{i} \lambda_i \sin\theta_i\ d\theta_i\land d\phi_i$$
and a Poisson structure given by
$$
\lbrace f,g \rbrace= \sum_{i} \frac{1}{\lambda_i}\ \left(\frac{\partial f}{\partial\cos\theta_i}\frac{\partial g}{\partial\phi_i}-\frac{\partial f}{\partial\phi_i}\frac{\partial g}{\partial\cos\theta_i}\right) .
$$
For any function $f$, we have $\dot f=\lbrace f,H \rbrace$.
Embedding $S^2$ in $\R^3$, the Poisson bracket becomes:
$$
\lbrace f,g \rbrace= -\sum_{i} \frac{1}{\lambda_i}\ \left(\frac{\partial f}{\partial x_i},\frac{\partial g}{\partial x_i},x_i \right)
$$
where $(a,b,c)=(a\times b)\cdot c$.

\subsection{Symmetries of the equations}
\label{symsection}
In this part, we apply the work of Lim, Montaldi, and Roberts \cite{LMR00}.
Let $\tau$ be the permutation of $S_{2N}$, $\tau=\prod_{i=1}^N \tau_{i,i^{\prime}}$ where $\tau_{i,j}$ is the transposition which exchanges $x_i$ and $x_j$.
That is, $\tau$ is a  permutation  of order two which exchanges ($+$)vortices  with ($-$)vortices.

Let ${G}=O(3)\times S_N \times S_N\rtimes \Z_2 \lbrack \tau \rbrack$ where $H\rtimes K$ is the semi-direct product of the groups $H$ and $K$.
The group ${G}$ acts on $\cal{P}$ in the following manner:
$$
g\cdot(x_1,\dots,x_{N},x_1^{\prime},\dots,x_{N}^{\prime})=(Ax_{\tau^k\sigma(1)},\dots,Ax_{\tau^k\sigma(N)},Ax_{\tau^k\sigma^{\prime}(1^{\prime})},\dots,Ax_{\tau^k\sigma^{\prime}(N^{\prime})})
$$
where $g=(A,\sigma,\sigma^\prime,\tau^k)\in {G}$.
Let $\chi$ be the following homomorphism,
\begin{eqnarray*}
\chi : {G} & \to & \lbrace -1,+1 \rbrace\\
(A,\sigma,\sigma^\prime,\tau^k) & \mapsto & (-1)^k\ \det A .
\end{eqnarray*}
The following proposition gives symmetries of the Hamiltonian and the dynamical system.
\begin{prop}
\ 

\noindent $\bullet$ The Hamiltonian $H$ is invariant under the action of $G$ on $\cal{P}$.\\
$\bullet$ Elements of $\Ker(\chi)$ are symmetries of $X_H$, while elements of $G\setminus \Ker(\chi)$ are time reversing symmetries of $X_H$.
\end{prop}
\proof
(i) Easy to check.\\
(ii) Let $g=(A,\sigma,\sigma^\prime,\tau^k)\in G$ and $f,h \in C^\infty(\mathcal{P})$, we have 
\begin{eqnarray*}
\lbrace f,h \rbrace (g \cdot x) & = & (-1)^{k+1} \sum_{i} \frac{1}{\lambda_i}\ \left(\frac{\partial f}{\partial x_i}(Ax),\frac{\partial h}{\partial x_i}(Ax),Ax_i \right)\\ 
 & = & (-1)^{k+1}\ \det A \sum_{i} \frac{1}{\lambda_i}\ \left(\frac{\partial fA}{\partial x_i},\frac{\partial hA}{\partial x_i},x_i\right) = \chi (g)\;\lbrace fg,hg \rbrace (x).
\end{eqnarray*}
Since $\lbrace f,H \rbrace (x)=df_x\cdot X_H(x)$, $$X_H(g\cdot x)=\chi(g)\; dg_x \cdot X_H(x).\Box$$

The action of $G$ is \textit{semi-symplectic} (see \cite{MoR} for the study of such an action in the reduced space): elements $g\in\Ker(\chi)$ are \textit{symplectic}, they preserve the symplectic form and the Poisson bracket, while  elements of $g\in G\setminus \Ker(\chi)$ are \textit{anti-symplectic}, the symplectic form and the Poisson bracket are transformed into their opposite.
The anti-symplectic elements are time reversing symmetries of $X_H$:
if $t\mapsto x(t)$ is a  solution of the dynamical system with $x(0)=x_0$, then the solution with initial value $g\cdot x_0$ is  $t\mapsto g\cdot x(-t)$. 

The isotropy subgroup ${G_x}=\lbrace g\in G \mid g\cdot x = x \rbrace$ of a point $x\in \cal{P}$ will be called  the \emph{symmetry group} of  $x$.
Let $K$ be a subgoup of $G$, we recall that the fixed point set of $K$ is $\Fix (K,\p)=\lbrace x\in\mathcal{P} \mid g\cdot x = x, \forall g\in K  \rbrace$, and that if $K$ is formed of symplectic elements, then $\Fix (K,\p)$ is invariant under the flow.
When there are no ambiguities, we write $\Fix (K)$ instead of $\Fix (K,\p)$.\\

\rmk
We can construct a symplectic element from two anti-symplectic elements (a reflection composed with the permutation $\tau$).
For example, consider $K$ the subgroup of $G$ generated by the two elements $(r_{\frac{2\pi}{N}},(1\dots N)^{-1},(1^{\prime}\dots N^{\prime})^{-1})$ and $(s_{x},\tau)$.
$\Fix (K)$ is the set of latitudinal semi-regular $2N$-gons formed of alternating vorticities.
The subgroup $K$ is formed of symplectic elements, $\Fix(K)$ is then invariant under the flow.
Actually in this case the dynamic can be completely described:
see \cite{To}.

\subsection{Conserved quantities}

From Emmy Noether's Theorem, Hamiltonian systems with continuous symmetry group satisfy conservation laws.
These conserved quantities are the components of the \emph{momentum map}.

Consider the diagonal action of $SO(3)$ on $\p$.
Let $\xi\in\so$ and $\xi_\p$ the vector field obtained by differentiating $g\cdot x$ with respect to $g$ in the direction $\xi$ at $g=e$.
Identifing $\so$ with $\R^3$ equipped with the cross product as Lie bracket, we get:
$$
\xi_\p(x)=(\xi\times x_1,\dots,\xi\times x_{2N})
$$
Recall that a \emph{momentum map} $J : \p\to \so^* \simeq\R^3$ for this action is such that
$$
\lbrace f,\left< J(x),\xi \right> \rbrace = df_x \cdot \xi_P 
$$
for all $f\in C^\infty(\p)$ and $\xi\in\so\simeq\R^3$.
We identify the natural bracket between $\so$ and its dual with the canonical scalar product of $\R^3$.

Let the \emph{center of vorticity vector} be the following vector:
$$
\Phi(x)=\sum_{j=1}^{2N}\lambda_j x_j .
$$
\begin{prop}
The momentum map for the action of SO(3) is $J(x)=-\Phi(x)$.
\end{prop}
\proof
It follows from the definition of momentum maps $J$ and the expression of $\xi_\p$ that for all $f\in C^\infty(\p)$,
$$
-\sum_{j} \frac{1}{\lambda_j}\ \left(\frac{\partial \<J(x),\xi\m}{\partial x_j}\times x_j \right) \cdot\frac{\partial f}{\partial x_j}=\sum_{j} \xi\times x_j \cdot\frac{\partial f}{\partial x_j} .
$$
Thus for all $\xi$ and $j$, $\lambda_j\xi\times x_j=-\frac{\partial \<J(x),\xi\m}{\partial x_j}\times x_j$, therefore $\frac{\partial J(x)}{\partial x_j}=-\lambda_j$ for all $j=1,\dots,2N$. 
Hence the momentum map is $J(x)=-\sum_j\lambda_j x_j$ up to a constant.
\ep

Consequently, the components of $\Phi$ are three conserved quantities (we could compute this directly differentiating $\Phi$ with respect to time).

Throughout this paper, we will take a frame $(O,\vec{e_x},\vec{e_y},\vec{e_z})$ of $\R^3$ such that $\Phi$ is parallel to the $(Oz)$ axis.
We obtain the following properties of $\Phi$ from its formula.
\begin{prop}
\ 

\noindent $\bullet$ $\Phi$ is $G$-equivariant with the following action on $\so^*\simeq\R^3$: 
$$g\cdot \zeta=(-1)^k\ A\zeta$$ where $g=(A,\sigma,\sigma^\prime,\tau^k)\in G$.

\noindent $\bullet$ Let $K$ be a subgroup of $G$.
If $x\in Fix(K,\p)$, then  $\Phi (x)\in  Fix(K,\R^3)$.
\end{prop}

\section{Equilibria}

\subsection{Fixed equilibria}

To emphasize the difference between equilibria and \riab, we call the former \textit{fixed equilibria}.
Thus a fixed equilibrium is a critical point of the Hamiltonian.
First, we give a general result which will serve to prove existence of fixed equilibria.
\begin{theo}
\label{fix}
Let $K$ be a subgroup of $G$.
An isolated point in $\Fix(K)$ is a fixed equilibrium.
\end{theo}
This Theorem was given first by Michel \cite{Mi71}, and it is also an application of the \textit{Principle of Symmetric Criticality} of Palais \cite{P79} which states that if the directional derivatives $dH_x(u)$ vanish for all directions $u$ at $x$ tangent to $\Fix(K)$, then directional derivatives in directions  transverse  to $\Fix(K)$ also vanish.
In particular, an isolated point of $\Fix(K)$ is a critical  point of $H$.

Note that equilibria derived by this Theorem do not depend on the form of the Hamiltonian, they depend only on the symmetry group and the phase space.
For example, the system of $N$ point charges on a sphere is Hamiltonian with potential (\cite{GE92}):
$$
W=\sum_{i<j} \lambda_i\lambda_j/l_{ij}
$$
where $\lambda_i$ is the charge of the $i$-th particle and $l_{ij}$ is the Euclidian distance between particles $i$ and $j$.
Hence, fixed equilibria of the vortex system described thereafter are also fixed equilibria of the point charges system.
But the stability results for fixed equilibria of the vortex system do not hold for the corresponding fixed equilibria of the point charges system, since the stability analysis depends on the expression for the Hamiltonian.

Kidambi and Newton \cite{KN98} found that a fixed equilibrium formed of vortices $x_i$ with vorticities $\lambda_i$ must statisfy $\sum_{i} (\lambda_i-\sum_j \lambda_j)\lambda_{i} x_{i} = 0$.
In the present case of $2N$ vortices $x_{1},\dots,x_{2N}$ with $\lambda_1=\dots=\lambda_N=-\lambda_{N+1}=\dots=-\lambda_{2N}$, this implies that \textit{the centre of the sphere is the isobarycentre of the configuration (barycentre with identical weights)}.
\begin{prop}
\label{fixexist}
The following  configurations are fixed equilibria:

\noindent - an equatorial ($\pm$)ring.

\noindent - a regular tetrahedron formed of $4$ ($+$)vortices and its tetrahedron dual formed of $4$ ($-$)vortices.
\end{prop}
Following the notation introduced in \cite{LMR00}, we denote the equatorial ($\pm$)ring by $D_{2Nh}(R_e)$ (see figure 1).\\
\begin{figure}[h]
    \begin{center}
\psfrag{+}{$+$}\psfrag{-}{$-$}
 \includegraphics[width=2.4in,height=1.933in,angle=0]{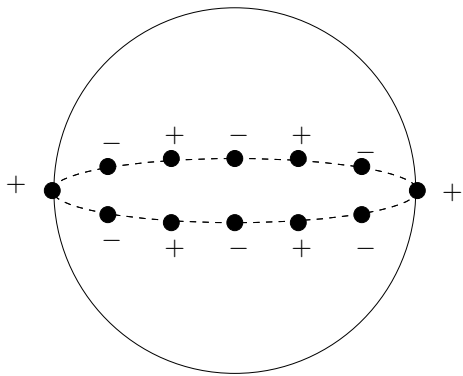}
 \caption{Equatorial ($\pm$)ring $D_{2Nh}(R_e)$.}
     \end{center}
 \end{figure}

\proof
Let $K$ be the symmetry goup of an equatorial ($\pm$)ring.
We have $\pi(K)=D_{2N} \times \Z_2[s_z]$ where $\pi$ is the cartesian projection $\pi:G \to O(3)$ and $D_{N}$ is the dihedral group of order $2N$.
Consider also $K^{\prime}=\left< (g,{\sigma}_{g}^{-1},{\sigma}_{g}^{-1}),g\in\T \right> \subset G$ where $\T$ is the subgroup of $O(3)$ formed of the symmetries of a regular tetrahedron, and  $\sigma_{g}$ is the permutation induced by $g$ on the vertices.
The reader can check that $\Fix(K)$ and $\Fix(K^{\prime})$ are formed of isolated points.
So $\Fix(K)$ and $\Fix(K^{\prime})$ are formed of fixed equilibria by use of Theorem \ref{fix}.
\ep

\rmk
This proposition is also true if the tetrahedron and its dual do not have opposite vorticities.
Indeed, the proof does not use the permutation $\tau$ which exchanges the $(+)$ with the ($-$).
In the same manner, the following arrangements are fixed equilibria: 
a cube formed of $8$ identical vortices and its octahedron dual formed of $6$ identical vortices, a dodecahedron formed of $20$ identical vortices and its icosahedron dual formed of $12$ identical vortices.
For a topological proof and some other fixed equilibria formed of regular polyhedra, see \cite{LMR00}.

\subsection{Relative equilibria}

Let $x_e\in\p$ and $x_e(t)$ the dynamic orbit of $X_H$ with $x_e(0)=x_e$.
Then $x_e$ is called a \textit{relative equilibrium} if for all $t$ there exists $g_t \in SO(3)$ such that $x_e(t)=g_t \cdot x_e$.

Let $\mu=\Phi(x_e) \in \R^3$.
We recall that $x_e$ is relative equilibrium if and only if there exists $\xi \in\so\simeq \R^3$ such that $x_e$ is a critical point of $$H_{\xi}(x)=H(x)+\left< \Phi(x)-\mu,\xi\right>$$ where the pairing $\left< \right>$ between $\R^3$ and its dual is the canonical scalar product of $\R^3$.
We can say equivalently that relative equilibria are critical points of the restriction of $H$ to $\Phi^{-1}(\mu)$ since the level set $\Phi^{-1}(\mu)$ are always non singular for point vortex systems of more than two vortices.

The solution is given by $x_e(t)=\exp(t\xi)\cdot x_e$, the vortices turn therefore uniformly around $\xi$ and $\xi$ is called the \textit{angular velocity} of the relative equilibrium $x_e$.

We can choose the action of $G$ on $\so$ such that if a relative equilibrum has a particular symmetry, then its angular velocity has also this symmetry, as the following result shows.
\begin{prop}
\label{symangvel}
Let $x_e$ be a relative equilibrium with angular velocity $\xi$.
Then $\xi$ satisfies $$\Ad_g\xi= \chi(g)\xi,\ \forall g\in G_{x_e}$$
where $\chi$ is defined in Section \ref{symsection}.
Moreover, $\xi\in \Fix(G_{x_e},\so)$ where the action of $G$ on $\so$ is defined by $g\cdot \xi = \chi(g)\Ad_g\xi $.
\end{prop}
\proof
Let $x$ be a relative equilibrium with angular velocity $\xi$, and $g\in G_{x}$.
Hence, $dH_x=\<\xi,d\Phi_x\m$.
Since $dH_x=dH_{g\cdot x}\cdot dg_x=dH_{x}\cdot dg_x$ and $d\Phi_{x}\cdot dg_x=\chi(g)\CoAd_g d\Phi_{x}$,
$$
\<\Ad_g\xi,d\Phi_x\m=\<\xi,\CoAd_{g^{-1}}d\Phi_{x}\m=\chi(g)\<\xi,d\Phi_x\cdot dg_x\m=\chi(g)\<\xi,d\Phi_x\m .
$$
Thus $\Ad_g\xi= \chi(g)\xi$ since $\Phi$ is a submersion.
Clearly $\xi\in \Fix(G_{x},\so)$ with respect to the action $g\cdot \xi = \chi(g)\Ad_g\xi $.
\ep

The following property is due to the rotational symmetry.
\begin{prop}
\label{angvel}
The angular velocity of a relative equilibrium is parallel to $\Phi$: the vortices rotate around the momentum vector.
\end{prop}
\proof
Let $x_e$ be a relative equilibrium with angular velocity $\xi$ and $\mu=\Phi(x_e)$.
We have $x_e(t)=\exp(t\xi)\cdot x_e\in \Phi^{-1}(\mu)\cap G\cdot x_e = G_\mu\cdot x_e$, so $\exp(t\xi)\in G_\mu$ and $\CoAd_{\exp(t\xi)}\mu=\mu$.
As $\CoAd_A=A$ for $A\in SO(3)$, we get $\exp(t\xi)\mu=\mu$.
A differentiation of the last identity will give the result. 
\ep

\rmk
As pointed out by Kidambi and Newton \cite{KN98}, ${\Phi}\cdot {x_i}$ depends only on the distances $l_{ij}$ for all $i$, and so these are conserved quantities if $x_1,\dots,x_{2N}$ is a \rumb.
Therefore each vortex moves on a cone around the momentum vector $\Phi$.\\

We recall that we defined $\vec{e_z}$ such that  $\vec{e_z}=\Phi/ \|\Phi\|$.
From Proposition \ref{angvel}, it follows that $\xi= \dot \phi_i \vec{e_z}$ for all $i=1,\dots,2N$, hence we obtain the following formula for the angular velocity.
\begin{prop}
\label{angvelexpr}
Let $x_e=(x_1,\dots,x_{2N})$ be a relative equilibrium.
The angular velocity $\xi$ satisfies:
$$
\xi=\frac{1}{\sin\theta_i} \sum_{j,j\neq i} \lambda_j \frac{\sin\theta_i \cos\theta_j -\sin\theta_j \cos\theta_i \cos(\phi_i-\phi_j)} {1-\cos\theta_i \cos\theta_j -\sin\theta_i \sin\theta_j \cos(\phi_i-\phi_j)}\  \vec{e_z}
$$
for all $i=1,\dots,2N$, where $\lambda_i$, $\theta_i$ and $\phi_i$ are respectively the vorticity, the co-latitude and  the longitude of $x_i$.
\end{prop}

\subsubsection{Large symmetry relative equilibria}

We first give a theorem which is the analogue of Theorem \ref{fix} for \riab. 
\begin{theo}
\label{releq}
Let $K$ be a subgroup of $G$, $x\in\Fix(K)\subset\p$ and $\mu=\Phi(x)$.
If $x$ is an isolated point in $\Fix(K)\cap\Phi^{-1}(\mu)$, then $x$ is a relative equilibrium.
If in addition $K$ is formed of symplectic elements, then $x$ is a fixed equilibrium.
\end{theo} 
The proof is an application of the \emph{Principle of Symmetric Criticality} \cite{P79} and the fact that a relative equilibrium is a critical point of the restriction of $H$ to a level set $\Phi^{-1}(\mu)$.
The last statement follows from the conservation of the momentum map, and the fact that $\Fix(K)$ is invariant under the flow if $K$ is symplectic.

\begin{define}
Let $x\in\mathcal{P}$, $K$ its symmetry group in $G$ and $\mu=\Phi(x)$.
A relative equilibrium $x$ is said to be a \emph{large symmetry relative equilibrium} if  $x$ is isolated in $\Fix(K)\cap \Phi^{-1}(\mu)$.
\end{define}

Large symmetry relative equilibria are relative equilibria with maximal isotropy subgroup at a fixed momentum value.
Note again that a large symmetry relative equilibrium depends only on the phase space, the momentum map and on the symmetries of the Hamiltonian: a large symmetry relative equilibrium is a critical point for all maps $f:\Phi^{-1}(\mu) \to \R$ which are $G$-invariant.

First, we focus on relative equilibria which lie in a great circle, since their dynamics should differ most from planar dynamics.
\begin{prop}
For $2N\geq 8$, every large symmetry relative equilibrium which lies in a great circle is a fixed equilibrium.
\end{prop}
\proof
To prove this proposition, we could classify all the possibilities by listing finite subgroups of $O(3)$ and examine the dimension of the fixed point set in each case.
Here is a proof using the properties of the angular velocity.

Fix a great circle, say $\lbrace y=0 \rbrace$.
Let $x$ be a relative equilibrium which lies in this great circle, and let $\xi$ be its angular velocity.
Its symmetry group $G_x$ can be written $\< K,s_y: y\mapsto -y \m$ where $K$ is a subgroup of $G$. 
Consider now the set of the ($+$)vortices $x_+=(x_1,\dots,x_N)$ and $\< K_+ ,s_y \m$ the symmetry group of $x_+$ where $K_+$ is a subgroup of $O(3)\times S_N$. 
Similarly, $\< K_-,s_y \m$ is the symmetry group of $x_-=(x_1^\prime,\dots,x_N^\prime)$.
Let $\pi$ be the canonical projection $\pi:K \to O(2)$.
Without loss of generality, assume that $\pi(K_+)\subset \pi(K_-)$.
Since $G_x=\< K, s_y \m$ and $K_+\subset K$,
$$
\xi\in \Fix(\< K, s_y \m,\R^3)\subset \Fix(\< K_+, s_y \m,\R^3)\ \ \mbox{(Prop. \ref{symangvel})}.
$$
Suppose $D_2\subset \pi(K_+)$, then $\Fix(\< K_+, s_y \m ,\R^3) \subset \Fix(\< D_2,s_y \m,\R^3)=0$.  
Now $\Fix(\< D_2,s_y \m,\R^3)$ is the intersection of three pairwise orthogonal planes, hence equal to $0$.
Thus $\Fix(\< K, s_y \m ,\R^3) \subset \Fix(\< K_+, s_y \m ,\R^3)=0$, $\xi$ is zero and every relative equilibrium is then a fixed equilibrium.

Cases $\pi(K_+)=0$ and $\pi(K_+)=\Z_2$ remain to be treated: in these two cases $\dim\Fix(G_x)\geq N-2$.
Assume now that $2N\geq 8$, then $\dim\Fix(G_x)\geq 2$.
Since $x$ is a great circle configuration, $\Fix(G_x)$ is a closed submanifold of $T^{2N}$ the $2N$-torus.
Hence $\dim\Fix(G_x)\cap \Phi^{-1}(\mu)=\dim\Fix(G_x)-1$ and then $\dim\Fix(G_x)\cap \Phi^{-1}(\mu)\geq 1$:
large symmetry relative equilibria can not exist.
\ep 

For $2N\leq 6$, we find easily that large symmetry relative equilibria which lie in a great circle are:
a (+)ring formed of 2 vortices and 2 ($-$)vortices at the poles  (denoted $C_{2v}(R,2p)$), two regular 2-rings of opposite vorticities and on opposite latitudes with possibly two polar vortices of opposite vorticities, and  the equatorial ($\pm$)rings formed of 4 and 6 vortices (the latters are fixed equilibria, see Prop. \ref{fixexist}).

The following proposition gives a necessary condition to have a relative equilibrium formed of two regular $N$-rings.
\begin{prop}[Lim-Montaldi-Roberts \cite{LMR00}]
Consider a relative equilibrium formed of two regular N-rings, each of identical vortices with non-zero vorticities, together with $k_p=0,1$ or $2$ poles, and such that $\xi\neq 0$ or $\mu\neq 0$.
Then:

\noindent - either the two rings are in phase

\noindent - or the two rings are out of phase with an offset equal to $\pi\over N$.

\noindent We denote these configurations respectively by $C_{Nv}(2R,k_pp)$ and $C_{Nv}(R,R',k_pp)$.
\end{prop}
\textbf{Explanation about notation}
Notations of relative equilibria are inspired by the Lim-Montaldi-Roberts orbit type notation: 
that is a relative equilibrium $x_e=(x_1,\dots,x_{2N})$ is denoted $\Gamma (k_e R_e, k R, k^\prime R', k_p p)$ where $\Gamma=\pi(G_{x_e})$ ($\pi$ is the cartesian projection $\pi:G \to O(3)$), $R_e$ is an equatorial ($\pm$)ring of $2N$ vortices, $R$ is a latitudinal $N$-ring of identical vortices, $R'$ is a latitudinal $N$-ring of identical vortices dual to $R$, $p$ is a polar vortex, and    $k_e$, $k$, $k^\prime$, $k_p$ number the occurence of $R_e$, $R$, $R'$, $p$ in the configuration $x_e$  (the \textit{Sch\"{o}nflies-Eyring} notation is used for subgroups of $O(3)$).

Now, we show that some of these configurations are relative equilibria.
\begin{theo}
\label{interm}
Configurations $C_{Nv}(2R,k_pp)$ and  $C_{Nv}(R,R',k_pp)$ with $k_p=0$ or $2$ poles of opposite vorticities and such that the rings have opposite latitudes and opposite vorticities, are relative equilibria.
We denote these configurations respectively by $D_{Nh}(2R,k_pp)$ and $D_{Nd}(R,R',k_pp)$.
\end{theo}
\begin{figure}[ht]
    \begin{center}
\psfrag{-}{$-$}
\psfrag{D_4h(2R)}{$D_{4h}(2R)$}
\psfrag{D_4d(R,R')}{$D_{4d}(R,R')$}
\includegraphics[width=1.7in,height=2.3in,angle=0]{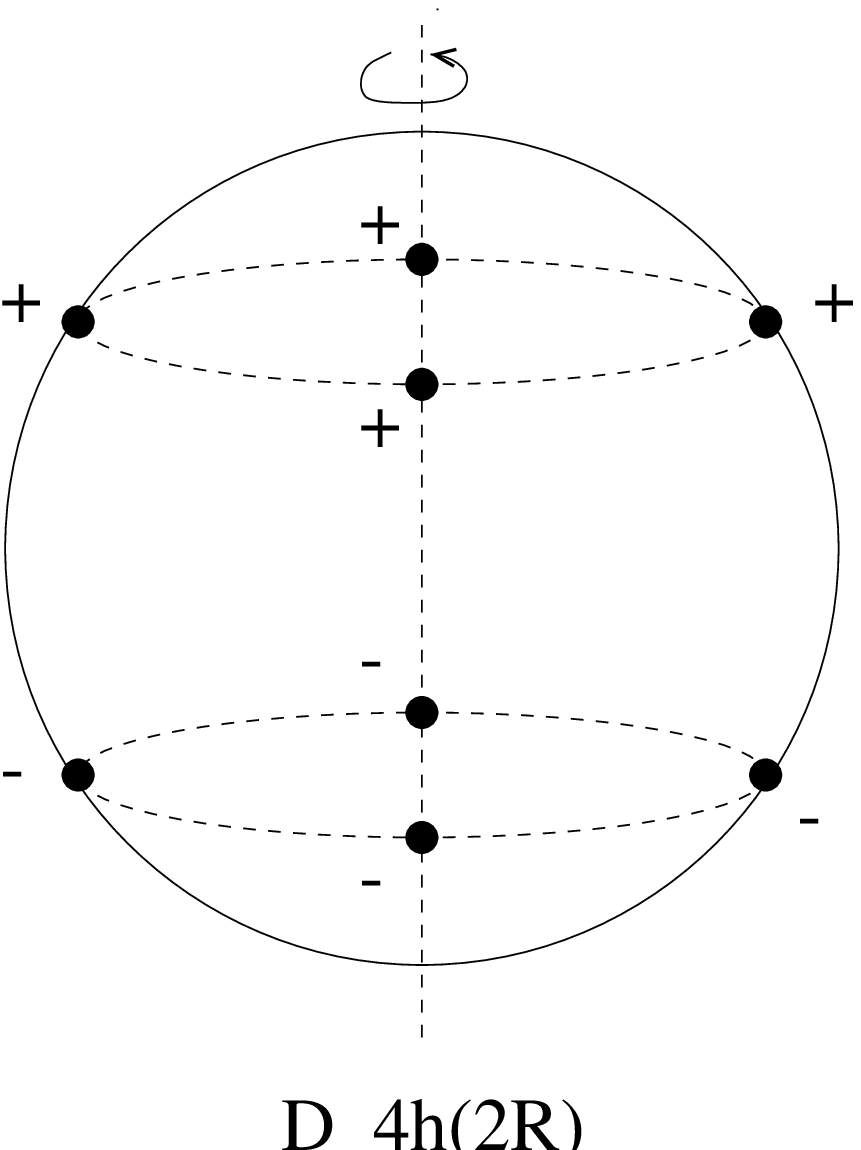}
\hspace{0.4in}
\includegraphics[width=1.7in,height=2.3in,angle=0]{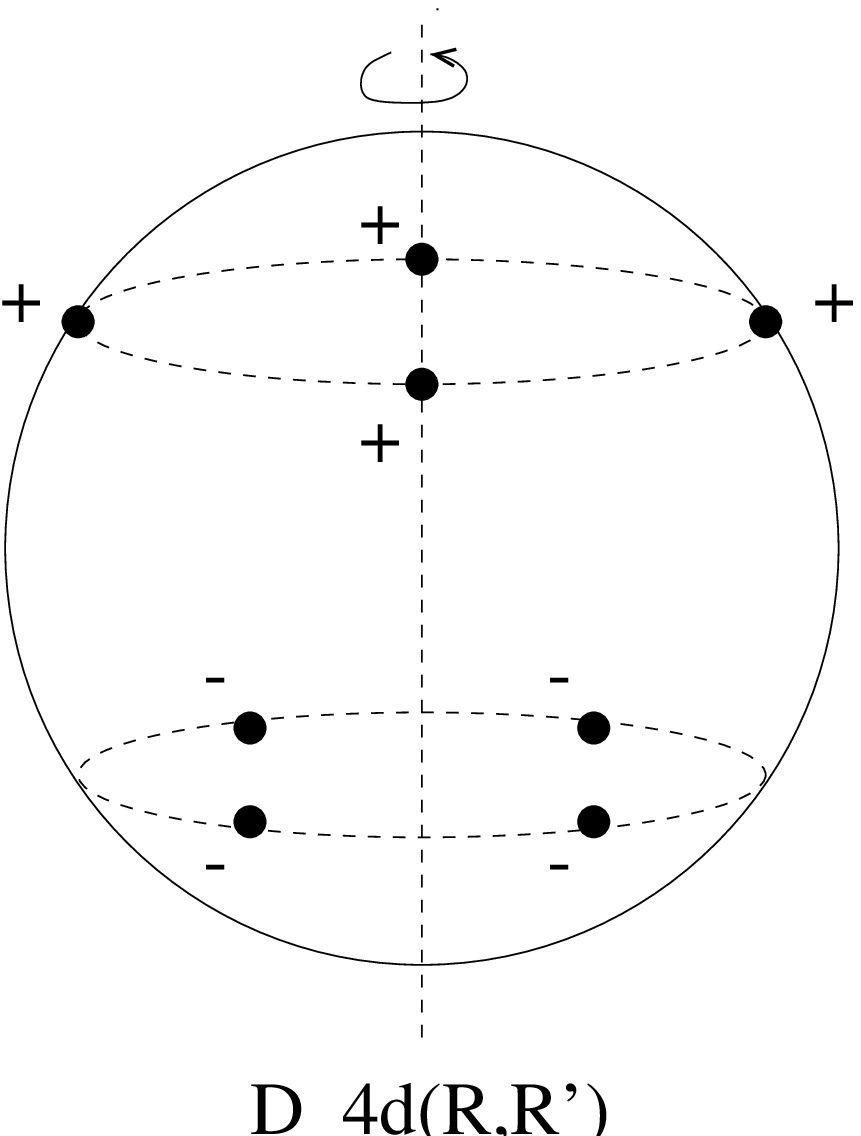}
\caption{Relative equilibria $D_{4h}(2R)$ and $D_{4d}(R,R')$.}
     \end{center}
 \end{figure}
\proof
For each configuration, the method is the following:
one computes $\Fix(K)\cap \Phi^{-1}(\mu)$ where $K$ is the isotropy subgroup of the configuration, and one uses Theorem \ref{releq}.
For example, consider a configuration $D_{Nd}(R,R')$ with $N$ odd.
Let $l$ be the integer part of $N/2$.
The isotropy subgroup $K$ of such a configuration is generated by the two elements $$(r_{\frac{\pi}{N}}\circ s_z,(1\dots N)^{-1},(1^\prime\dots N^\prime)^{-1},\tau ),\ (s_x,\tau_{2,N}\cdots\tau_{l+1,l+2},\tau_{1^\prime,N^\prime}\cdots\tau_{l^\prime,(l+2)^\prime}).$$
The submanifold $\Fix(K)$ is parametrized by $\theta_1$ and hence $\Fix(K)\cap\Phi^{-1}(\mu)$ is a singleton.
\ep

\rmk
Let $N$ be even and $K$ be the symmetry group of two aligned semi-regular $N$-gons of opposite vorticities and on opposite latitudes.
$K$ is isomorphic to $D_{Nh}$, $\Fix(K)\simeq\{(\theta_1,\phi_2)\in S^1\times S^1\}$ and $\Phi_{\mid \Fix(K)}(\theta_1,\phi_2)= 2N\cos\theta_1$ (make a picture!).
For given $\mu$, $-H$ restricted to $\Fix(K)\cap\Phi^{-1}(\mu)$ goes to infinity when $\phi_2$ goes to $\phi_1$ or $\phi_3$ which are fixed.
It follows that $H$ must have a maximum for a certain value of $\phi_2$.
In fact, $H$ is maximal for $\phi_2=(\phi_1-\phi_3)/2$ and this corresponds to the $D_{Nh}(2R)$ configuration.
This idea can be generalized and leads to \textit{minimal strata} \cite{LMR00}.\\

The theorem above does not assume that the polar and ring vortex strengths coincide, and in fact it is interesting to be able to vary their ratio.
We therefore fix the ring vorticities to be $\pm 1$ as usual, and let $\lambda_n$ and  $\lambda_s$ be the strengths of the North pole and South pole respectively, with $\lambda_n=-\lambda_s$.
The following proposition gives formulae for the angular velocity of \ria $D_{Nh}(2R,k_pp)$ and $D_{Nd}(R,R',k_pp)$.
\begin{prop}
The angular velocity of $D_{Nh}(2R,k_pp)$ relative equilibria is:
$$
{\xi}=\cos\theta\ \left[ \frac{N-1}{\sin^2\theta} + \sum_{j=1}^{N} \frac{1 + \cos\frac{2\pi}{N}(j-1)} {2-\sin^2\theta\ (1 + \cos\frac{2\pi}{N}(j-1))}+\frac{k_p\lambda_n}{\cos\theta \sin^2\theta}\right]\ \vec{e_z}.
$$
The angular velocity of $D_{Nd}(R,R',k_pp)$ relative equilibria is:
$$
{\xi}=\cos\theta\ \left[ \frac{N-1}{\sin^2\theta} + \sum_{j=1}^{N} \frac{1 + \cos(\frac{2\pi}{N}(j-1)+\frac{\pi}{N})} {2-\sin^2\theta\ (1 + \cos(\frac{2\pi}{N}(j-1)+\frac{\pi}{N}))}+\frac{k_p\lambda_n}{\cos\theta \sin^2\theta}\right]\ \vec{e_z}.
$$
In both cases $\theta$ is the co-latitude of the ($+$)ring.
\end{prop}
\begin{figure}[ht]
    \begin{center}
\psfrag{pi}{$\frac{\pi}{2}$}
\psfrag{l0}{$\lambda_n=0$}
\psfrag{lm}{$\lambda_n=-1$}
\psfrag{lp}{$\lambda_n=+1$}
\includegraphics[width=1.5in,height=2in,angle=0]{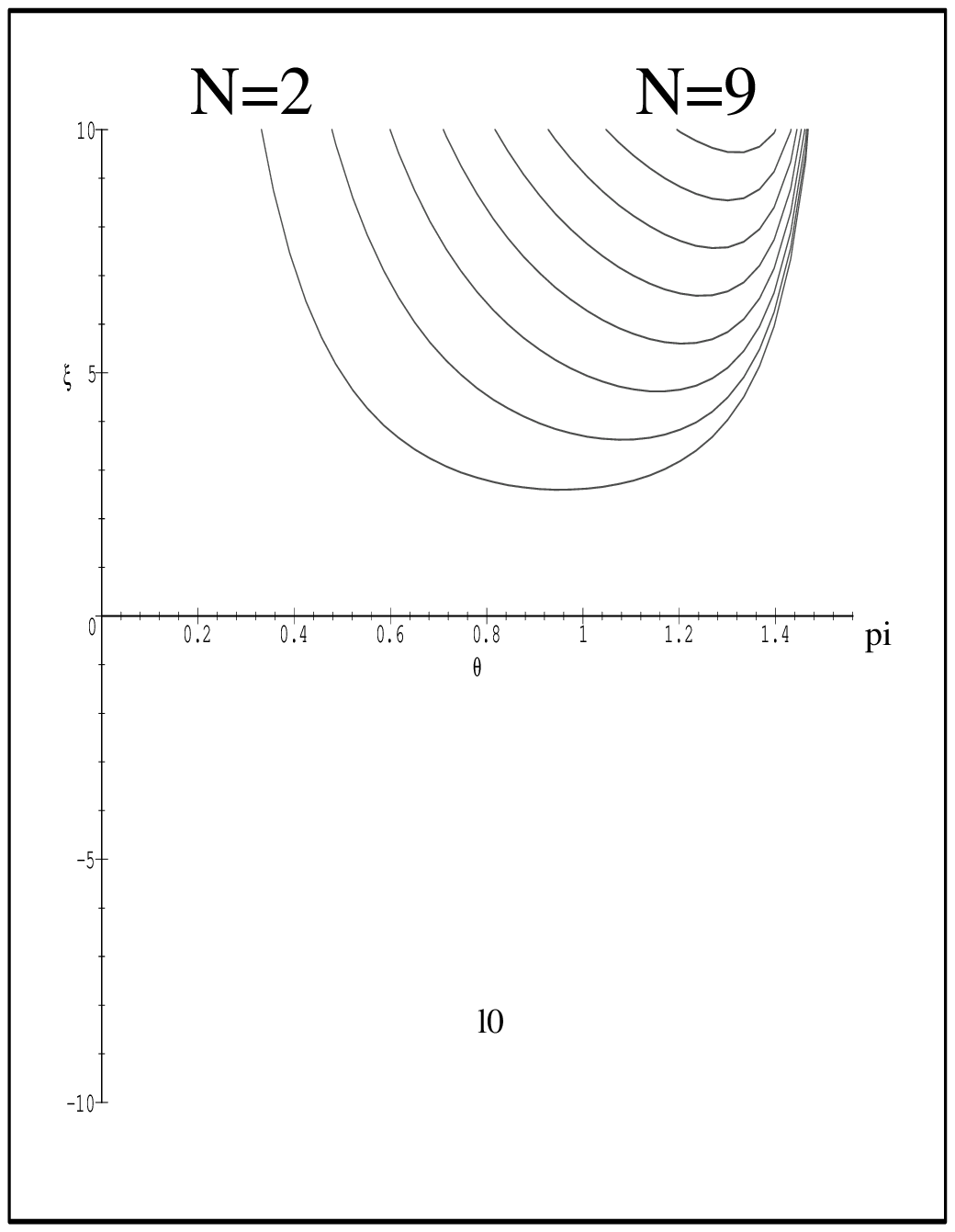}
\includegraphics[width=1.5in,height=2in,angle=0]{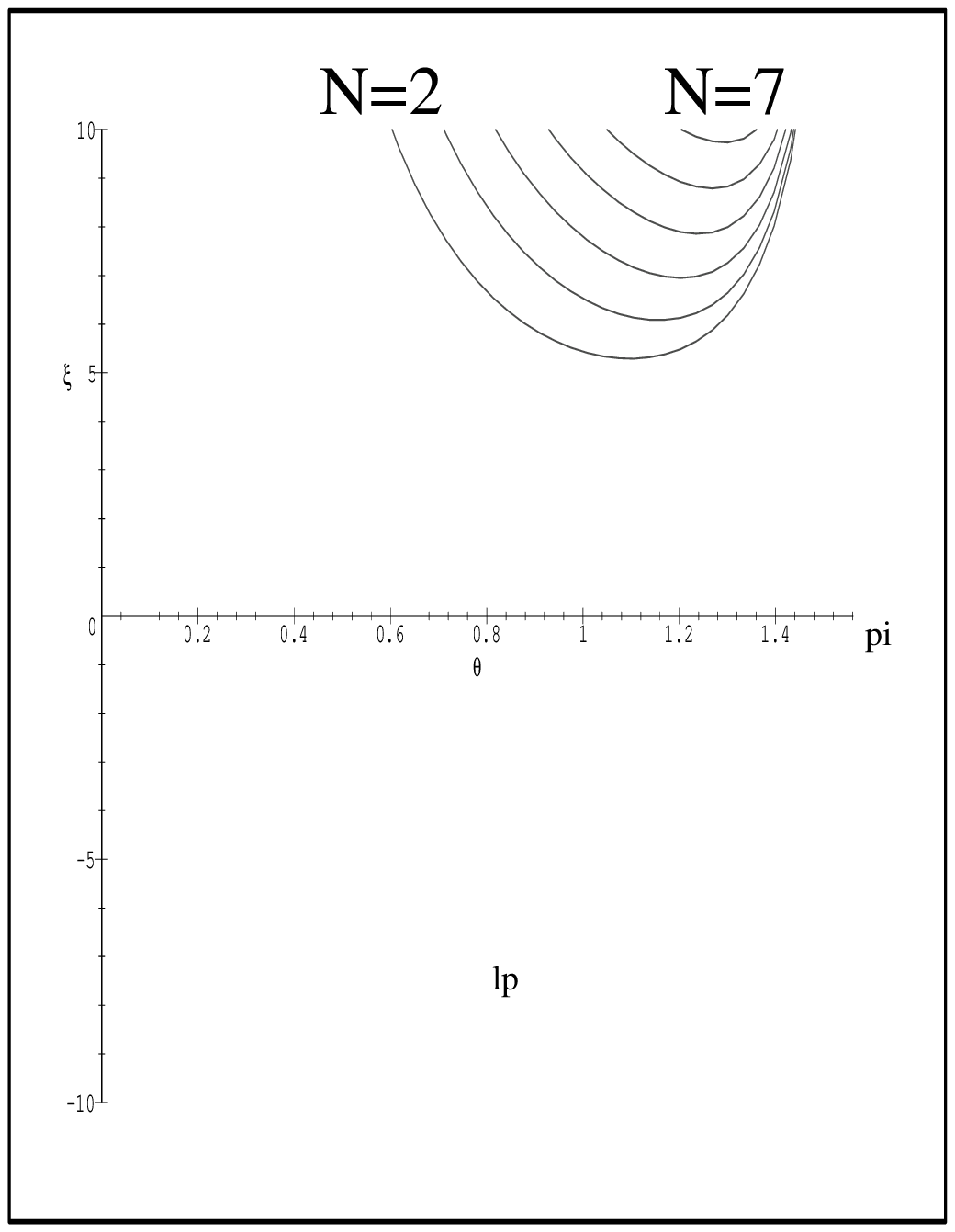}
\includegraphics[width=1.5in,height=2in,angle=0]{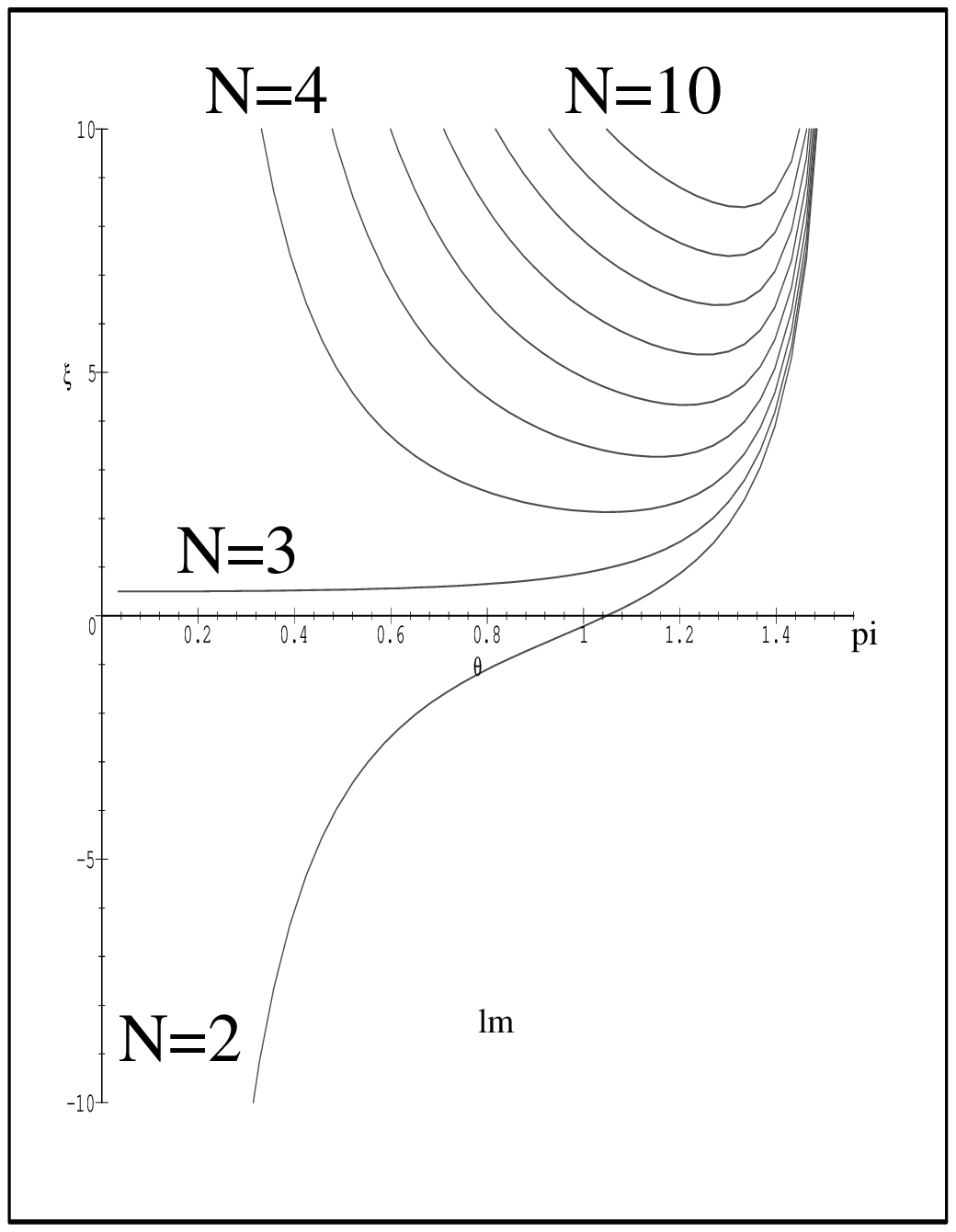}     
     \caption{Angular velocity of $D_{Nh}(2R,k_pp)$.}
     \end{center}
 \end{figure}
\begin{figure}[ht]
    \begin{center}
\psfrag{pi}{$\frac{\pi}{2}$}
\psfrag{l0}{$\lambda_n=0$}
\psfrag{lm}{$\lambda_n=-1$}
\psfrag{lp}{$\lambda_n=+1$}
\includegraphics[width=1.5in,height=2in,angle=0]{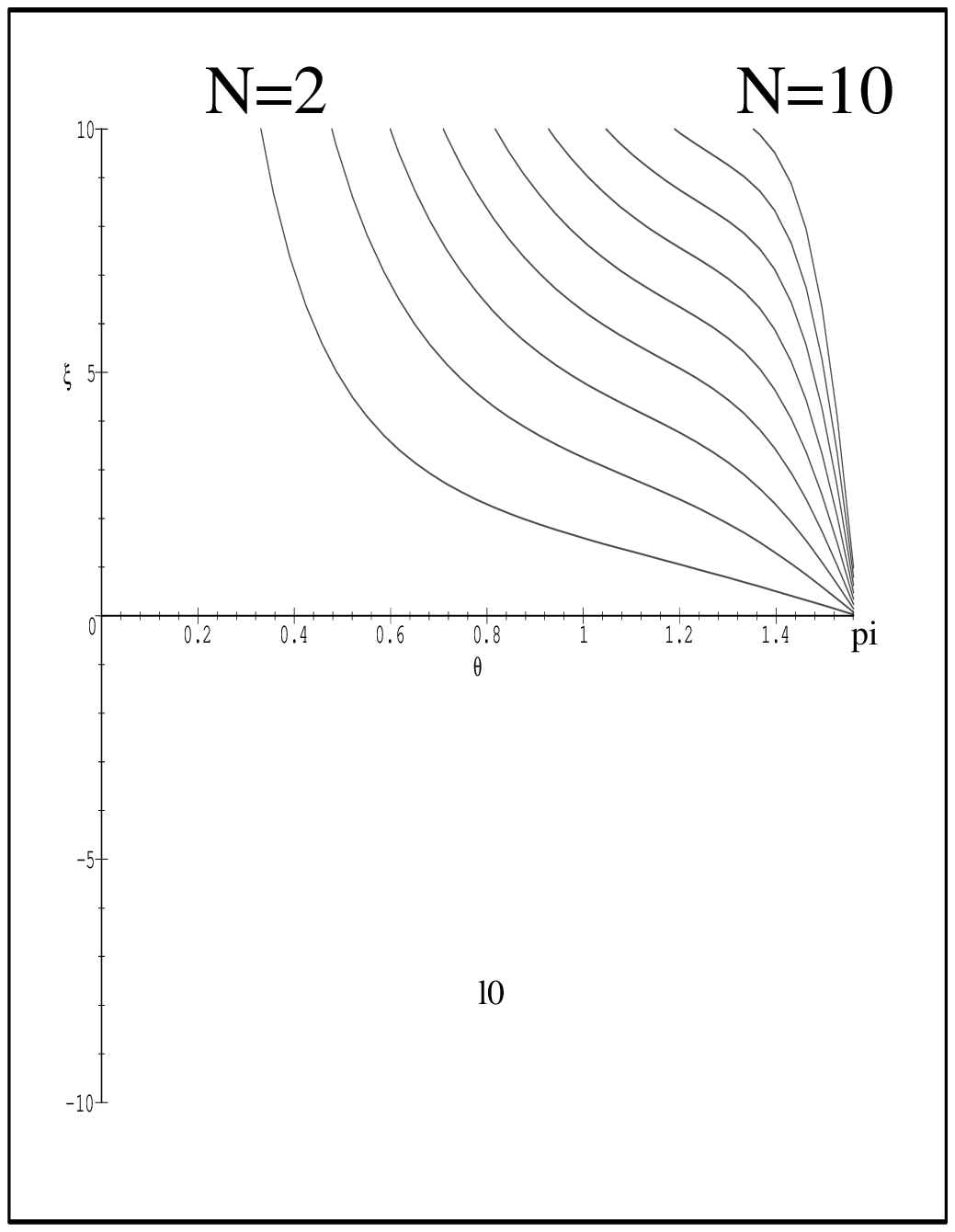}
\includegraphics[width=1.5in,height=2in,angle=0]{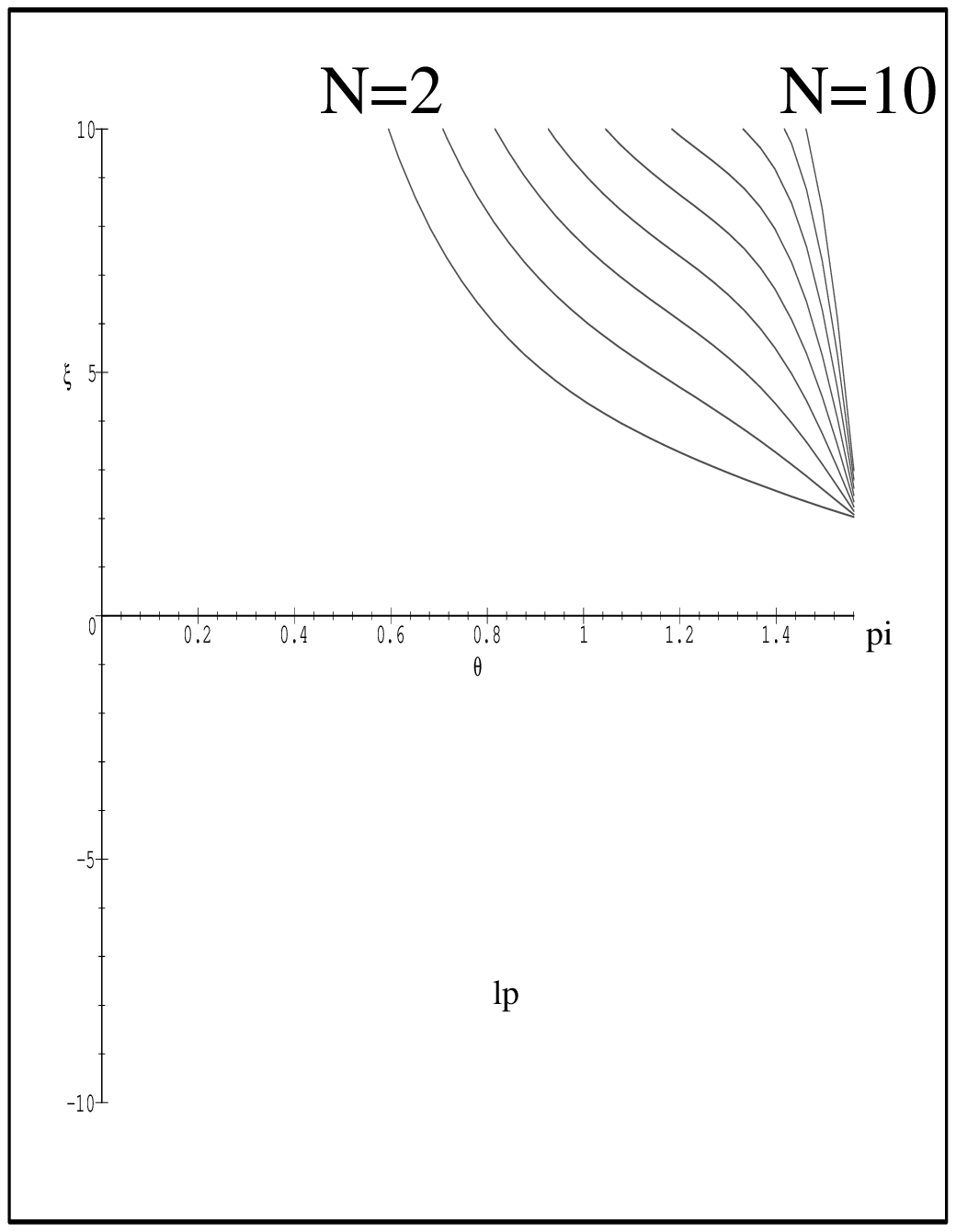}
\includegraphics[width=1.5in,height=2in,angle=0]{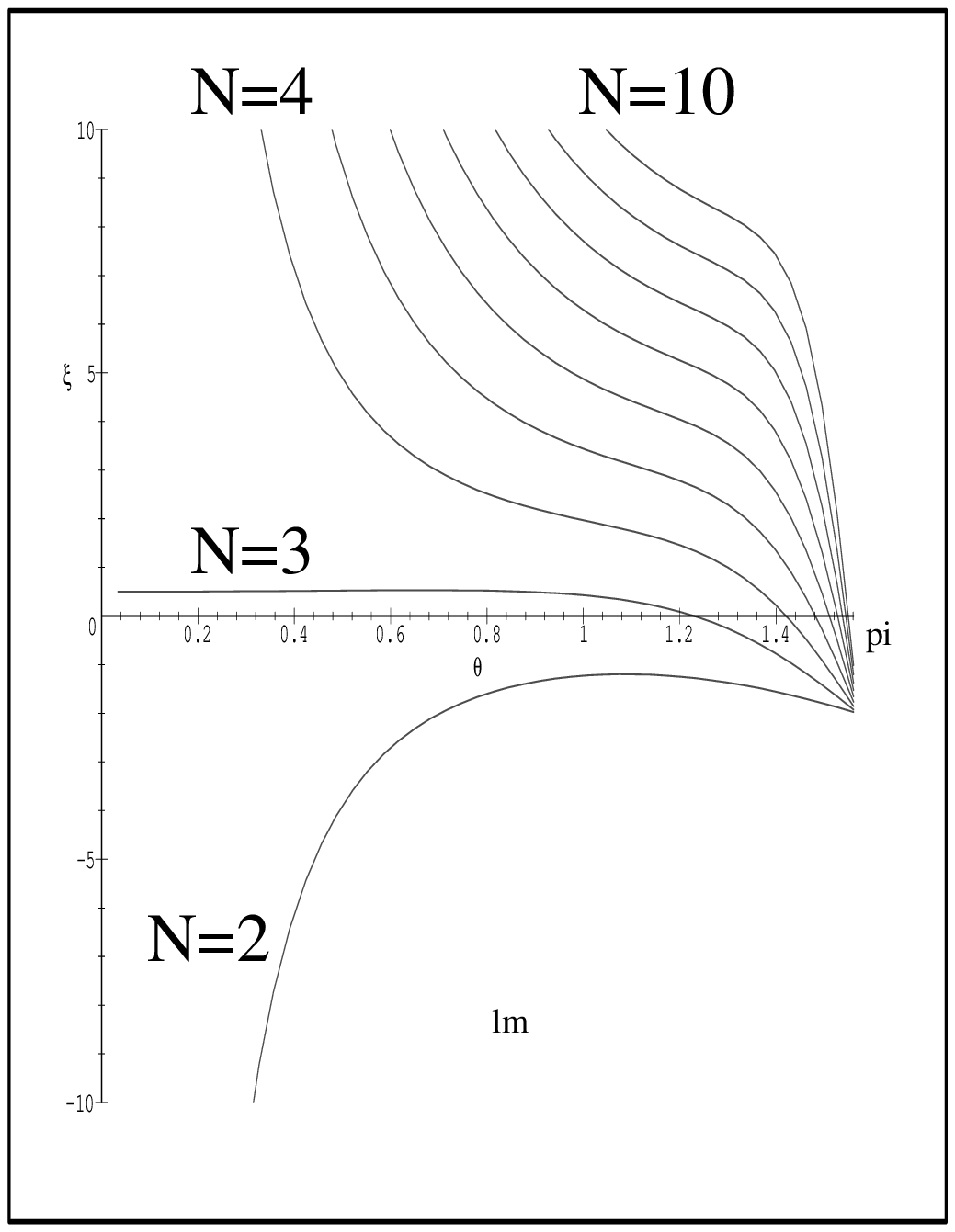}
\caption{Angular velocity of $D_{Nd}(2R,k_pp)$.}
     \end{center}
 \end{figure}
\proof
Let $x_e$ be a $D_{Nh}(2R,k_pp)$ or $D_{Nd}(2R,k_pp)$ \rumb, and $\theta$ be the co-latitude of the ($+$)ring in $x_e$.
We get from Prop. \ref{angvelexpr} that:
$$
\xi=\frac{1}{\sin\theta_1} \sum_{j=2}^{2N+k_p} \lambda_j \frac{\sin\theta_1 \cos\theta_j -\sin\theta_j \cos\theta_1 \cos(\phi_1-\phi_j)} {1-\cos\theta_1 \cos\theta_j -\sin\theta_1 \sin\theta_j \cos(\phi_1-\phi_j)}\  \vec{e_z}.
$$
Since $\theta_j=\tet$, $\tet_{N+j}=\pi-\tet$, $\lambda_j=+1$ and $\lambda_{N+j}=-1$ for all $j=1,\dots,N$, this expression becomes:
$$
\xi=\cos\tet\ \left[ \sum_{j=2}^{N} \frac{1-\cos\phi_j} {\sin^2\tet (1-\cos\phi_j)} - \sum_{j=1}^{N} \frac{-1-\cos\phi_{N+j}} {1+\cos^2\tet+\sin^2\tet \cos\phi_{N+j}}\right] + \frac{k_p\lambda_n}{\sin^2\theta}\ \vec{e_z}
$$
which leads to
$$
\xi=\cos\tet\ \left[ \frac{N-1}{\sin^2\theta} + \sum_{j=1}^{N} \frac{1 + \cos\phi_{N+j}}{2-\sin^2\theta\ (1 + \cos\phi_{N+j})} + \frac{k_p\lambda_n}{\cos\tet\sin^2\theta}\right]\ \vec{e_z}.
$$
The formula for $\xi$ follows from $\phi_{N+j}=\frac{2\pi}{N}(j-1)$ (resp. $\frac{2\pi}{N}(j-1)+\frac{\pi}{N}$) for a $D_{Nh}(2R,k_pp)$ (resp. $D_{Nd}(2R,k_pp)$) \rumb. 
\ep

For $\lambda_n\geq 0$, the angular velocity of $D_{Nh}$ goes to infinity as the two rings approach the  equator or poles, while the angular velocity of $D_{Nd}$ goes to infinity as the two rings approach the poles and goes to $k_p\lambda_n$ as the two rings approach the equator, we find again that an equatorial ($\pm$)ring is a fixed equilibrium.
For $\lambda_n <0$, other behaviour appears for small $N$ (see figures 3 and 4 for $\lambda_n=-1$).
Note that $\cos\theta\frac{N-1}{\sin^2\theta}$ is the angular velocity of relative equilibria formed of a single ring of $N$ identical vortices \cite{PD93,LMR}.

\subsubsection{Relative equilibria with less symmetry}
\label{lesssym}

In this section, we give relative equilibria with $\dim\Fix(K)=2$ for systems of 4 and 6 vortices, in order to find bifurcating branches of relative equilibria. 

First, we consider configurations of $6$ vortices formed of a (+)ring $R_+$ and a ($-$)ring $R_-$ together with two poles with vorticities $\lambda_n$, $\lambda_s$ such that $\lambda_s=-\lambda_n$.
That is,  we consider points of $\p$ with symmetry group $K=\left< (r_\pi,\tau_{1,2},\tau_{3,4})\right>$ where ``1,2'' label the (+)vortices, ``3,4'' the ($-$)vortices, and ``5,6'' the poles.
Let $\theta_+$ be the co-latitude of $R_+$ and $\theta_-$ that of $R_-$.  
The reader can check that we can not apply Theorem \ref{releq} since $\Fix(K)\cap\Phi^{-1}(\mu)$ has no isolated points.
However a critical point of $H$ in $\Fix(K)\cap\Phi^{-1}(\mu)$ is a relative equilibrium by the \textit{Principle of Symmetric Criticality}.
We therefore use the Lagrange multipliers Theorem to find critical points of $H_{\mid\Fix(K)}$ on the manifold $\Phi_{\mid\Fix(K)}^{-1}(\mu)$ in order to determine \riab.

For the configurations in question, the Hamiltonian on $\Fix(K)$ is deduced from (\ref{HAM}):
$$
H=\ln \left( \frac{(1-x^2)(1-y^2)}{ ((1-xy)^2-(1-x^2)(1-y^2)\cos^2\alpha )^2} \left[\frac{(1-x)(1+y)}{2(1+x)(1-y)}\right] ^{2\lambda_n} \right)
$$
where $x=\cos\theta_+$ and  $y=\cos\theta_-$ and $\alpha$ is the offset between $R$ and $R'$ defined modulo $\pi$. 
The momentum map on $\Fix(K)$ is given by $\Phi=2(x-y+\lambda_n)$.
Applying the Lagrange multipliers Theorem we find that there are two families of \ria with isotropy containing $K$:

\noindent - The first family satisfies $\alpha=0$ and 
$$
(y+x)[\; xy-1+2\lambda_n(y-x)\;]=0 . 
$$
The solution $y=-x$ correspond to the $D_{2h}(2R,2p)$ relative equilibria.
For $\vert \lambda_n \vert \geq \frac{1}{2}$, we have a branch of relative equilibria   $y=(2\lambda_n x+1)/(x+2\lambda_n)$ denoted $C_{2v}(2R,2p)$ and invariant under the symmetry $(x,y)\mapsto(-y,-x)$ (see figure 5). 
\begin{figure}[h]
    \begin{center}
\psfrag{ln=+1}{$\lambda_n=+1$}
\psfrag{ln=-0.6}{$\lambda_n=-0.6$}
\psfrag{D2h(2R,2p)}{$D_{2h}(2R,2p)$}
\psfrag{C2v(2R,2p)}{$C_{2v}(2R,2p)$}
\includegraphics[width=2.3in,height=2.3in,angle=0]{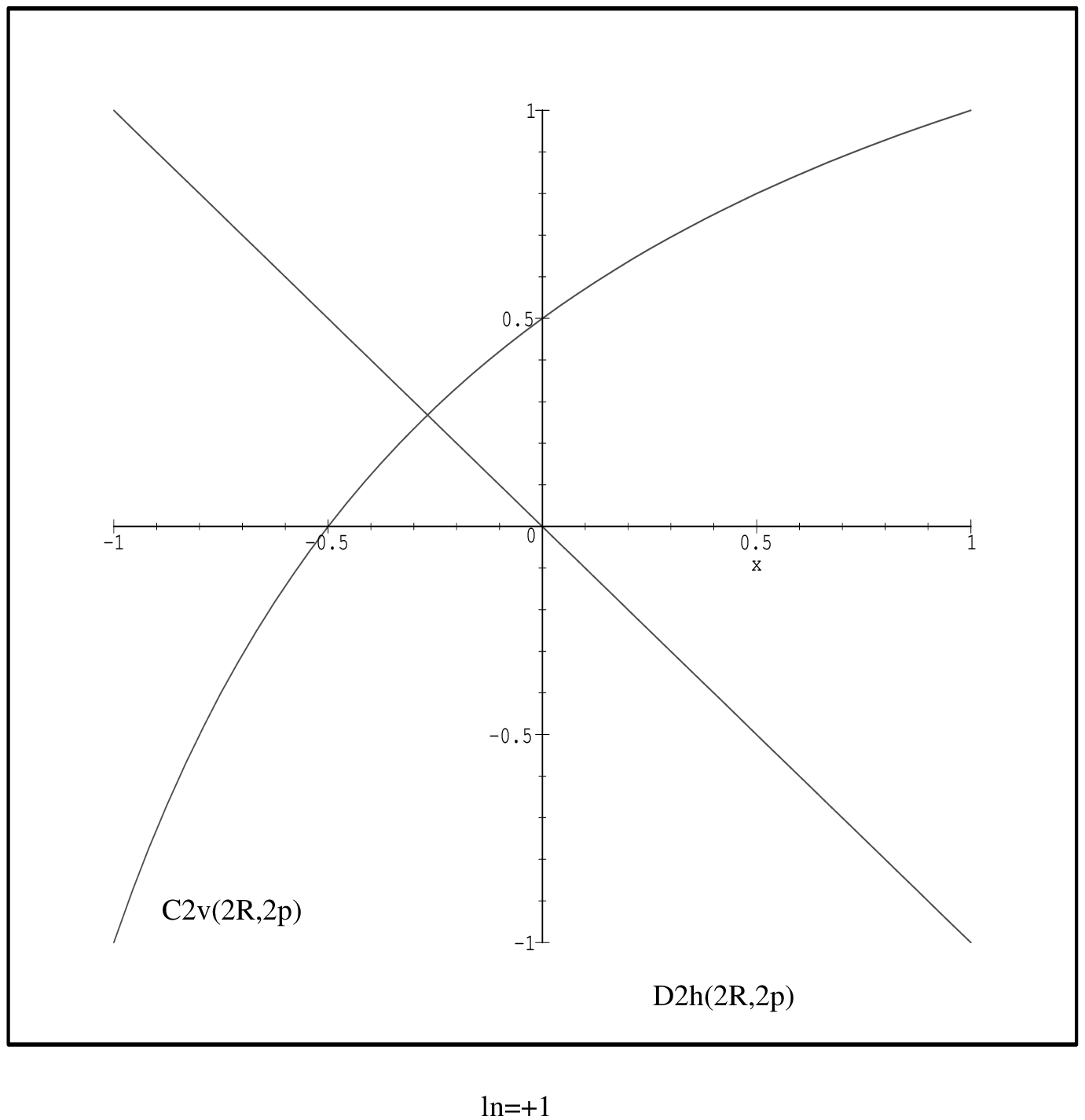}
\psfrag{D2h(2R,2p)}{$D_{2h}(2R,2p)$}
\psfrag{C2v(2R,2p)}{$C_{2v}(2R,2p)$}
\includegraphics[width=2.3in,height=2.3in,angle=0]{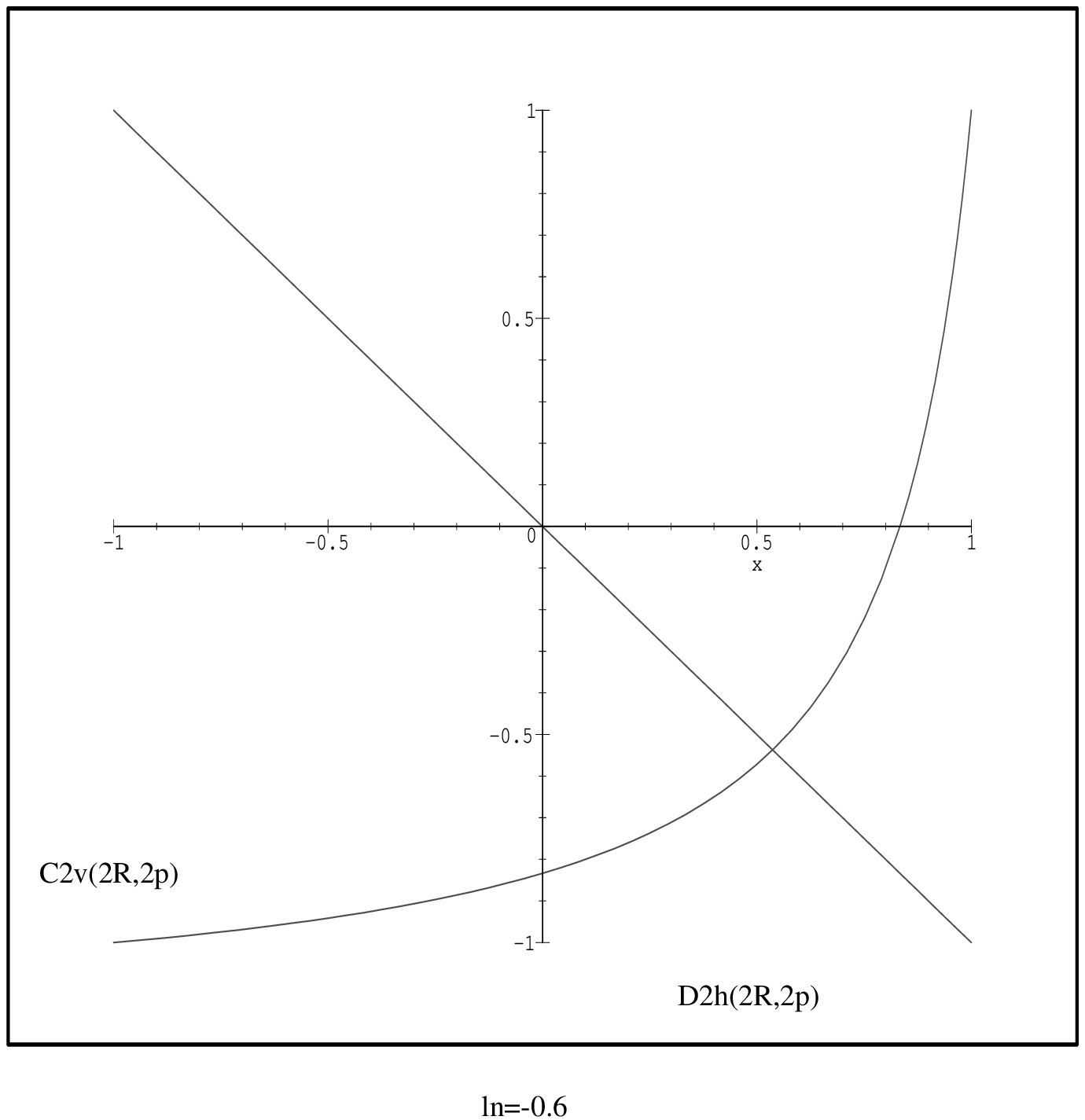}
\caption{Relative equilibria $D_{2h}(2R,2p)$ and $C_{2v}(2R,2p)$.}
     \end{center}
 \end{figure}
\begin{figure}[h]
    \begin{center}
\psfrag{ln=0}{$\lambda_n=0$}
\psfrag{ln=0.6}{$\lambda_n=0.6$}
\psfrag{D2d(R,R')}{$D_{2d}(R,R')$}
\psfrag{C2v(R,R')}{$C_{2v}(R,R')$}
\includegraphics[width=2.3in,height=2.3in,angle=0]{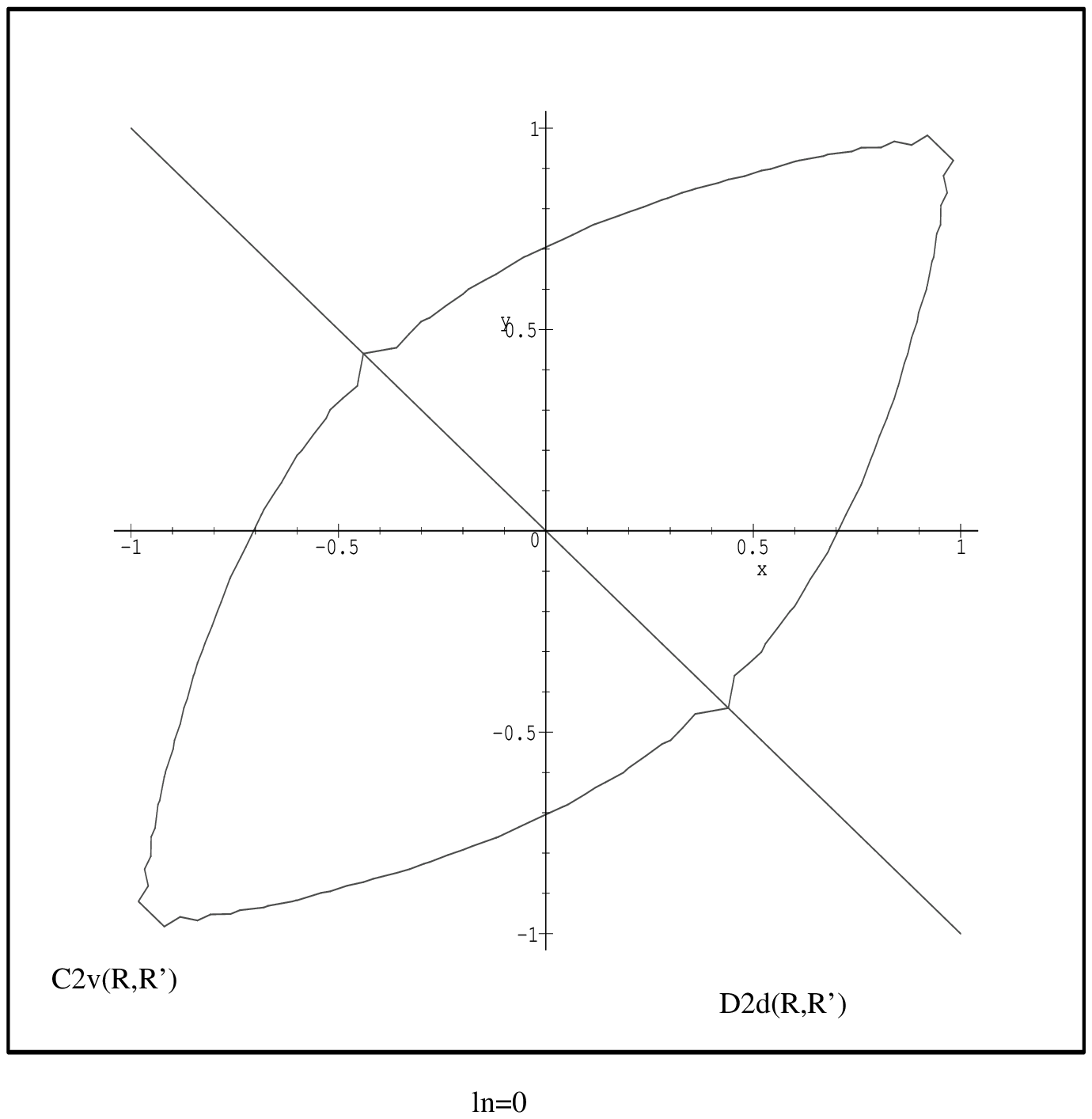}
\psfrag{D2d(R,R',2p)}{$D_{2d}(R,R',2p)$}
\psfrag{C2v(R,R',2p)}{$C_{2v}(R,R',2p)$}
\includegraphics[width=2.3in,height=2.3in,angle=0]{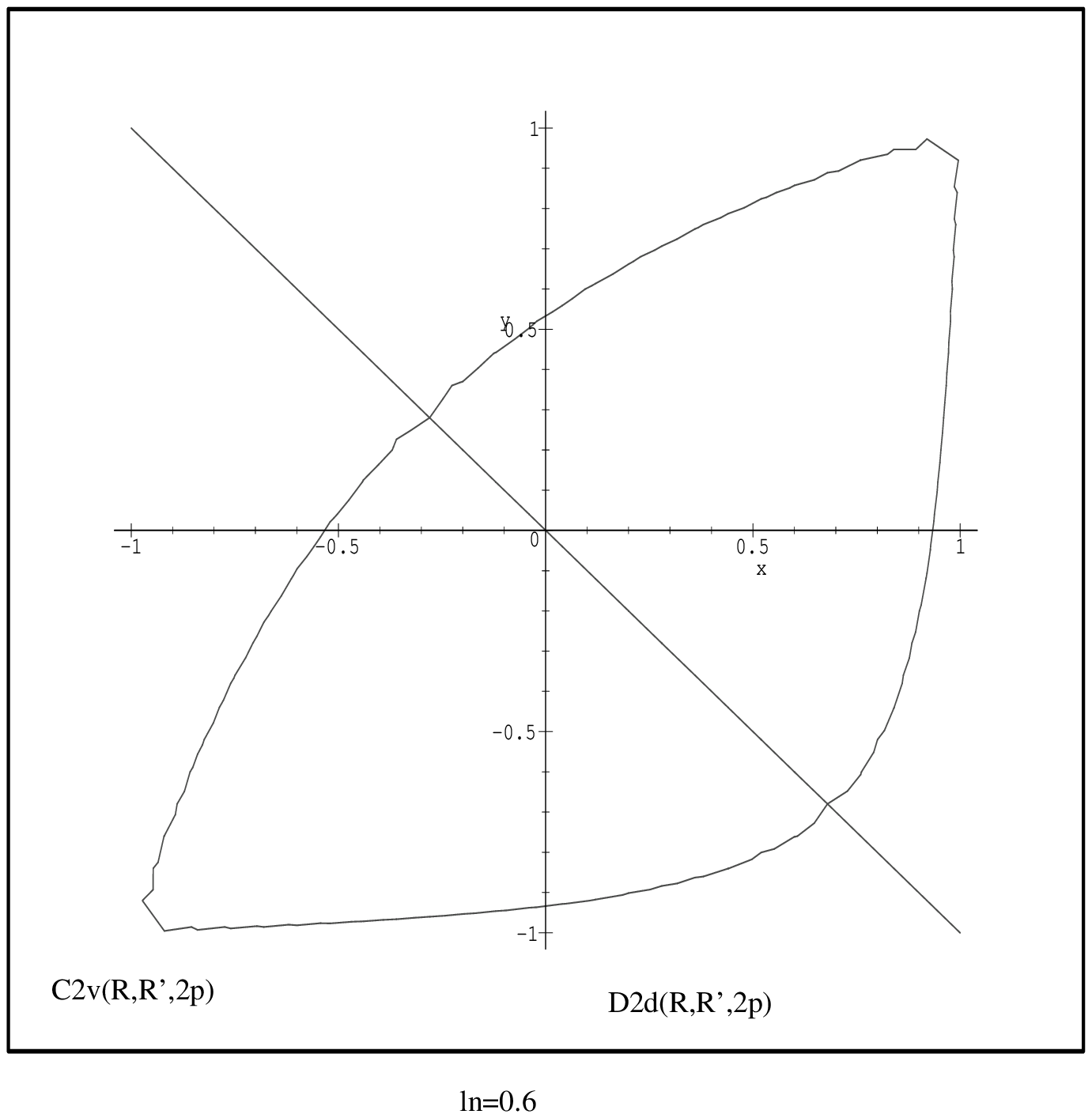}
\caption{Relative equilibria $D_{2d}(R,R',2p)$ and $C_{2v}(R,R',2p)$.}
     \end{center}
 \end{figure}

\noindent - The second family satifies $\alpha=\frac{\pi}{2}$ and 
$$
(y+x)\left[\; x^2y^2-2y^2-2x^2+2xy+1-2\lambda_n(1-xy)(x-y)\;\right]=0 .
$$
The branch $y=-x$ corresponds to the $D_{2d}(R,R^\prime,2p)$ relative equilibria.
The quartic curve  intersects the square $[-1,1]^2$ for all $\lambda_n$ so provides relative equilibria for any value of $\lambda_n$.
That is $$y=-\frac{x+\lambda_n(x^2+1)\pm (1-x^2)\sqrt{2+\lambda_n^2}}{x^2-2(1+\lambda_nx)}$$ is a branch of relative equilibria denoted $C_{2v}(R,R^\prime,2p)$. Note that $C_{2v}(R,R^\prime,2p)$ is invariant under $(x,y)\mapsto(-y,-x)$ and when $\lambda_n=0$ under one more symmetry $(x,y)\mapsto(y,x)$.
When $\lambda_n$ goes to infinity, $C_{2v}(R,R^\prime,2p)$ goes to the solution $y=x$ that is a  ($\pm$)ring with two poles of infinite and opposite strength.
See figure 6.\\

Second, we look for relative equilibria formed of two (+) and two ($-$) vortices with symmetry group $K=\left< (s_z:z\mapsto -z,\tau_{1,4}\tau_{2,3})\right>$ where ``1,2'' label the (+)vortices  and ``3,4'' the ($-$)ones.
So we have two longitudinal ($\pm$)rings $R_m,R_m^\prime$. 
We proceed as in the previous case: the Hamiltonian on $\Fix(K)$ is given by:
$$
H=2\ln\left( \frac{1-\sqrt{(1-x^2)(1-y^2)}\cos\alpha+xy}{4xy(1-\sqrt{(1-x^2)(1-y^2)}\cos\alpha-xy)}\right)
$$
where $x=\cos\theta_1$ and  $y=\cos\theta_3$ and $\alpha=\phi_3-\phi_1$ is the offset between $R_m$ and $R_m^\prime$ defined modulo $2\pi$. 
The momentum map on $\Fix(K)$ is given by $\Phi=2(x-y)$.
We find that there are two families of \ria with isotropy containing $K$:

\noindent - The first family satisfies $\alpha=0$ and  
$$
(y+x)[\; 2(yx^3+xy^3-x^2-y^2-xy+1) + (x^2+y^2+2xy-2)\sqrt{(1-x^2)(1-y^2)}\; ]=0 .
$$
But this does not provide relative equilibria since we have vortex collapse with $y=-x$ and the second term  does not vanish for $(x,y)\in [-1,1]^2$.

\noindent - The second family satisfies $\alpha=\pi $ and 
$$
(y+x)[\; 2(yx^3+xy^3-x^2-y^2-xy+1) - (x^2+y^2+2xy-2)\sqrt{(1-x^2)(1-y^2)}\; ]=0 .
$$
The branch $y=-x$ corresponds to the $D_{2h}(2R)$ relative equilibria.
The second term corresponds to a branch of relative equilibria denoted $C_{2v}(R_m,R_m^\prime)$ which is invariant under symmetries $(x,y)\mapsto(-y,-x)$ and  $(x,y)\mapsto(y,x)$.
When $x=y$ this branch coincides with  the $D_{4h}(R_e)$ fixed equilibria.
See figure 7.\\
\begin{figure}[h]
    \begin{center}
\psfrag{D4h(Re)}{$D_{4h}(R_e)$}
\psfrag{D2h(2R)}{$D_{2h}(2R)$}
\psfrag{C2v(Rm,Rm')}{$C_{2v}(R_m,R_m^\prime)$}
\includegraphics[width=3in,height=3in,angle=0]{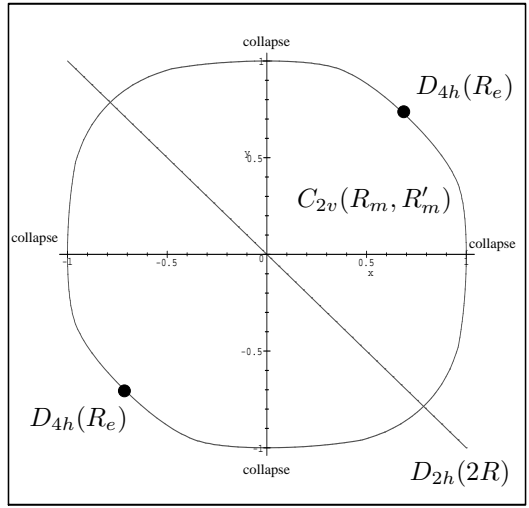}
\caption{Branches of relative equilibria $D_{2h}(2R)$ and $C_{2v}(R_m,R_m^\prime)$.}
     \end{center}
 \end{figure}

\rmks

\noindent $\bullet$ The equations for relative equilibria $x_e$ with $\dim G_{x_e}=n$ represent the intersection of $(n-1)$ hypersurfaces of $\R^n$.

\noindent $\bullet$ One can also prove existence of relative equilibria using works of Montaldi and Roberts \cite{Mo97,MoR99}:
let $x_e$ be a non-degenerate relative equilibrium (that is a non-degenerate critical point) such that $\Phi(x_e)=0$, and $x$ an arrangement with momentum $\mu$ near $0$ such that $G_x\subset G_{x_e}$.
Let $\cal{O}_\mu$ be the corresponding coadjoint orbit, that is $\cal{O}_\mu$ is the sphere of centre $O$ and radius $\vert\mu\vert$.
One has the following:
if $x$ is an isolated point of $\Fix(G_x,\cal{O}_\mu)$, then $x$ is a relative equilibrium.
It is an easy exercise to use this result in order to show that there exists a relative equilibrium of type $C_{2v}(R_m,R_m^\prime)$ near the fixed equilibrium $D_{4h}(R_e)$ (see figure 7).

\section{Stability of equilibria}

In order to determine stability of equilibria, we use the \textbf{energy momentum method} (see Marsden \cite{Ma92} and references therein for a general description).
Let $x_e$ be a relative  equilibrium, $\mu=\Phi(x_e)\in \R^3$ and $\xi$ its angular velocity  (if $x_e$ is a fixed equilibrium, $\xi=0$). The method consists of determining the symplectic slice $\G= T_{x_e}(SO(3)_{\mu}\cdot x_e)^\bot \cap \Ker D\Phi(x_e)$ (so $\G$ is transversal to $SO(3)_{\mu}\cdot x_e$) where
$$
SO(3)_{\mu}=\lbrace g\in SO(3)\mid\CoAd_g\cdot \mu=\mu \rbrace
$$
and then examine wether $d^2H_{\xi}|_{\G}(x_e)$ is definite or not where
$$
H_{\xi}(x)=H(x)+\<\Phi(x)-\mu,\xi\m.
$$
Note that if $d^2H_{\xi}|_{\G}(x_e)$ is  definite, that is $d^2H_{\xi}|_{\G}(x_e)$ is  positive definite or  negative definite, then $x_e$ is an extremum of the restriction of $H_{\xi}$ to the symplectic slice.

Let $K$ be a subgroup of $G$.
Recall that a relative  equilibrium $x_e$ is said to be \emph{Lyapunov stable modulo K}, if for all $K$-invariant open neighborhoods $V$ of $K\cdot x_e$ there is an open neighborhood $U\subseteq V$ of $x_e$ which is invariant under the Hamiltonian evolution.
Papers of Patrick \cite{Pa92} and more recently Lerman and Singer \cite{LS98} (see also \cite{OR99}) refined the \textbf{energy momentum theorem} given in \cite{Ma92}.
The assumptions of Patrick's Theorem are satisfied since $SO(3)$ is compact, we have then the following:

\noindent\textit{If $d^2H_{\xi}|_{\G}(x_e)$ is definite, then $x_e$ is Lyapunov stable modulo $SO(3)_{\mu}$.}\\
\\
For $\mu\neq 0$, $SO(3)_{\mu}$ is the set of rotations with axis $\<\mu\m$ and then isomorphic to $SO(2)$, while $SO(3)_{\mu} = SO(3)$  for $\mu = 0$.
Lyapunov stability modulo $SO(3)_{\mu}$ of a \rum with non zero angular velocity corresponds to the ordinary stability of the corresponding periodic orbit if $SO(3)_{\mu} = SO(2)$, that is if $\mu\neq 0$.
In the remaining of the section, we will omit \textit{modulo $SO(3)_{\mu}$}.

\subsection{Stability of  $D_{Nd}(R,R',k_pp)$ and $D_{Nh}(2R,k_pp)$ ($k_p=0,2$)}

The most general results of this section are \textit{linear instability} results: relative equilibria of type $D_{Nd}(R,R')$, $D_{Nh}(2R)$, and $D_{Nh}(2R,2p)$ are linearly unstable provided $N\geq 4$, $N\geq 7$, and $N\geq 9$ respectively.
The important fact is that linear instability implies nonlinear instability.
For $N\leq 6$ (resp. $N\leq 8$), \ria of type $D_{Nh}(2R)$ (resp. $D_{Nh}(2R,2p)$) are Lyapunov stable for a certain range of co-latitude.
First, one needs to describe the symplectic slice.

\textbf{Notation} 
Let $x_e$ be a relative equilibrium of type $D_{Nd}(R,R',k_pp)$ or $D_{Nh}(2R,k_pp)$ with $k_p=0,2$, and $\theta_0$ be the co-latitude of the (+)ring of $x_e$.   

Let $j=1,\dots,N$ number the vortices in cyclic order in the two rings.
Write tangent vectors $(\delta\theta_j,\delta\phi_j)$ at the $j^{th}$ (+)vortex, $(\delta\theta_j^\prime,\delta\phi_j^\prime)$ at the $j^{th}$ ($-$)vortex, $(\delta x_n,\delta y_n)$ at the North pole, and $(\delta x_s,\delta y_s)$ at the South pole.

Let $q$ be an integer, $0\leq q\leq N-1$, and $i=\sqrt{-1}$.
Define the following vectors of $T_{x_e}\p\oplus iT_{x_e}\p$:

$$\alpha_{q,\theta}+i\beta_{q,\theta}=\sum_{j=1}^{N} e^{\frac{2i\pi q}{N}(j-1)}(\delta\theta_j+e^{iq\varphi_0}\delta\theta_j^\prime) $$
$$\alpha_{q,\theta}^{\prime}+i\beta_{q,\theta}^{\prime}=\sum_{j=1}^{N} e^{\frac{2i\pi q}{N}(j-1)}(\delta\theta_j-e^{iq\varphi_0}\delta\theta_j^\prime) $$
$$\alpha_{q,\phi}+i\beta_{q,\phi}=\sum_{j=1}^{N} e^{\frac{2i\pi q}{N}(j-1)}(\delta\phi_j+e^{iq\varphi_0}\delta\phi_j^\prime) $$
$$\alpha_{q,\phi}^{\prime}+i\beta_{q,\phi}^{\prime}=\sum_{j=1}^{N} e^{\frac{2i\pi q}{N}(j-1)}(\delta\phi_j-e^{iq\varphi_0}\delta\phi_j^\prime) $$
where $\varphi_0$ is the offset between the two rings, that is $\varphi_0=0$ (resp. $\frac{\pi}{N}$) for a $D_{Nh}$ (resp. $D_{Nd}$) relative equilibrium.
Note that $\beta_{q,\theta}$, $\beta_{q,\theta^{\prime}}$, $\beta_{q,\phi}$ and $\beta_{q,\phi^{\prime}}$ vanish for $q=0$ and $N/2$ (for $N$ even) while others do not vanish.

Define then $\delta x_p=\delta x_n+\delta x_s$, $\delta y_p=\delta y_n+\delta y_s$, $\delta x_p^{\prime}=\delta x_n-\delta x_s$, and $\delta y_p^{\prime}=\delta y_n-\delta y_s$.

\subsubsection{Description of the symplectic slice}
\label{sympslice}

The symplectic slice $\G$ is $T_{x_e}(SO(3)_{\mu}\cdot x_e)^\bot \cap \Ker D\Phi(x_e)$.
In order to compute a basis of $\G$, one needs expressions of the differential of the momentum map and the tangent space to the group orbit at the relative equilibria.
\begin{prop}
\label{diffmom}
Let $x_e$ be a $D_{Nd}(R,R',k_pp)$ or $D_{Nh}(2R,k_pp)$ relative equilibrium.
The differential of the momentum map at $x_e$ is given by:
$$
D\Phi(x_e)(\delta\theta_j,\delta\theta_j^{\prime},\delta\phi_j,\delta\phi_j^{\prime},\delta x_n,\delta x_s,\delta y_n,\delta y_s)=$$
$$\cos(\theta_0)(\alpha_{1,\theta}+i\beta_{1,\theta})+i\sin(\theta_0)(\alpha_{1,\phi}^{\prime}+i\beta_{1,\phi}^{\prime})+\lambda_n(\delta x_p^{\prime}+i\delta y_p^{\prime})\oplus -\sin(\theta_0)\alpha_{0,\theta}^{\prime}
$$ 
where the direct sum corresponds to the $C_{nv}$-invariant decomposition of $\so^*$ as a direct sum of a plane and the line $\Fix(C_{nv},\so^*)$.
\end{prop}
\proof
The expression of the momentum map is 
$$\Phi(\theta_j,\theta_j^{\prime},\phi_j,\phi_j^{\prime},x_n,x_s,y_n,y_s)=\Phi_{r}+\Phi_{p}$$
where $$\begin{array}{l}
\Phi_{r} = \sum_{j=1}^{N}(\sin\theta_j\cos\phi_j-\sin\theta_j^{\prime}\cos\phi_j^{\prime},\sin\theta_j\sin\phi_j-\sin\theta_j^{\prime}\sin\phi_j^{\prime},\cos\theta_j-\cos\theta_j^{\prime})\\
\Phi_{p} = \lambda_n(x_n-x_s,y_n-y_s,\sqrt{x_n^2+y_n^2}+\sqrt{x_s^2+y_s^2}).
\end{array}$$
Then identifying $\so^*$ with $\C\oplus \< e_z \m$, we get $\Phi=\Phi_{\C}\oplus\Phi_{z}$
where $$\begin{array}{lll}
\Phi_{\C} & = & \sum_{j=1}^{N}(\sin\theta_j e^{i\phi_j}-\sin\theta_j^{\prime} e^{i\phi_j^{\prime}})+\lambda_n(x_n-x_s+i(y_n-y_s))\\
\Phi_z & = & \sum_{j=1}^{N} \cos\theta_j-\cos\theta_j^{\prime} + \lambda_n(\sqrt{x_n^2+y_n^2}+\sqrt{x_s^2+y_s^2}).
\end{array}$$
Differentiation of this expression at a $D_{Nd}(R,R',k_pp)$ or $D_{Nh}(2R,k_pp)$ configuration gives the result.
\ep
\begin{prop}
\label{tgtspace}
Let $x_e$ be a $D_{Nd}(R,R',k_pp)$ or $D_{Nh}(2R,k_pp)$ configuration, and $\mu=\Phi(x_e)$.\\
If $\mu\neq 0$, then $\so_{\mu} \cdot x_e$ is generated by the vector $\alpha_{0,\phi}$.\\ 
If $\mu = 0$, then $\so_{\mu} \cdot x_e$ is generated by the three vectors:
$$
\alpha_{0,\phi},\ \beta_{1,\theta}+\cos\theta_0\sin\theta_0\;\alpha_{1,\phi}^{\prime}+\delta y_p^{\prime},\ \alpha_{1,\theta}-\cos\theta_0\sin\theta_0\;\beta_{1,\phi}^{\prime}+\delta x_p^{\prime}
.$$
\end{prop}
\proof
Let $x_e$ be a $D_{Nd}(R,R',k_pp)$ or $D_{Nh}(2R,k_pp)$ arrangement.
The tangent space to the group orbit at $x_e$ is:
$$
\so\cdot x_e=\{(\eta\times x_1,\dots,\eta\times x_{2N+k_p})\mid \eta\in\R^3\} 
$$
Hence $\so\cdot x_e$ is generated by the three following vectors of $T_{x_e}\PP$:
$$
(e_x\times x_1,\dots,e_x\times x_{2N+k_p}),(e_y\times x_1,\dots,e_y\times x_{2N+k_p}),(e_z\times x_1,\dots,e_z\times x_{2N+k_p}).
$$
If $\mu=\Phi(x_e)\neq 0$, $SO(3)_\mu$ is the set of rotations about the $z$-axis  then $\so_{\mu} \cdot x_e$ is generated only by the last one $(e_z\times x_1,\dots,e_z\times x_{2N+k_p})$.
The expressions of these vectors in spherical coordinates give the desired results.
\ep 

If $\mu\neq 0$, it follows that the dimension of the symplectic slice $\G$ is $4N-4$ where $N$ is the number of (+)vortices, while if $\mu=0$ it is of dimension $4N-6$.

\subsubsection{Isotypic decomposition of the symplectic slice}

This second step consists of obtaining a block diagonalization of $d^2H_{\xi}|_{\G}(x_e)$ in order to compute more easily eigenvalues and then examine definiteness.
This is done by using the following method of group representation theory (see also \cite{GSS88} and \cite{CL00}).  

Let $V$ be a vector space and $G$ be a compact Lie group.
Recall that a representation $W\subset V$ of $G$ is said to be an \textit{irreducible representation} if $W$  has no proper $G$-invariant subspaces.
There are a finite number of distinct irreducible representations of $G$ in $V$, say $U_1,\dots,U_l$.
Let $V_k$ be the sum of all irreducible representations $W\subset V$ such that $W$ and $U_k$ are isomorphic representations.
Then $V=V_1\oplus\cdots\oplus V_l$.
This decomposition of $V$ is unique and is called \textit{the isotypic decomposition} \cite{Se78}.

Irreducible representations of groups $D_{Nd}$ and $D_{Nh}$ are either of dimension 1 or 2.
The following two propositions give irreducible representations of dimension 1 and 2 respectively.
\begin{prop}
\label{rep1}
Let $x_e$ be a $D_{Nd}(R,R',k_pp)$ or $D_{Nh}(2R,k_pp)$ configuration and $k_p=0,2$. 
There are four irreducible representations of dimension $1$ of $G_{x_e}$ on $T_{x_e}\p$:
$$\<\alpha_{0,\theta}\m_{\R} ,\; \<\alpha_{0,\phi}\m_{\R} ,\; \<\alpha_{0,\theta}^{\prime} \m_{\R} ,\; \<\alpha_{0,\phi}^{\prime} \m_{\R} \; .$$
Of these, only $\<\alpha_{0,\theta}\m_{\R}$ and $\<\alpha_{0,\phi}^{\prime}\m_{\R}$ lie in the symplectic slice.
\end{prop}
\proof
We give only the proof for $D_{Nd}(R,R',k_pp)$ arrangements as the one for $D_{Nh}(2R,k_pp)$ is quite similar.
The group $D_{Nd}$ is generated by two elements $a,b$ such that $a^2=1$, $b^{2N}=1$.
Let $w\in T_{x_e}\PP$ such that $a\cdot w=\lambda w$ and $b\cdot w=\mu w$ with $\lambda,\mu\in\R$. 
Hence $\lambda^2=1$, $\mu^{2N}=1$ and $\lambda=\pm1$, $\mu=\pm1$.
For $1\leq i\leq N$, $b^2\cdot(\delta\theta_i,\delta\phi_i)=(\delta\theta_{\sigma(i)},\delta\phi_{\sigma(i)})$ where $\sigma$ is a $N$-cycle of $S_N$.
Thus $(\delta\theta_{\sigma(i)},\delta\phi_{\sigma(i)})=\mu^2(\delta\theta_i,\delta\phi_i)=(\delta\theta_i,\delta\phi_i)$.
It follows that $\delta\theta_i=\delta\theta_j$ and $\delta\phi_i=\delta\phi_j$ for all $i,j=1,\dots,N$.
The same argument holds also for $\delta\theta_i^{\prime},\delta\phi_i^{\prime}$.
We have $b\cdot (\delta\theta_1,\delta\phi_1)=(-\delta\theta_1^{\prime},\delta\phi_1^{\prime})=\mu(\delta\theta_1,\delta\phi_1)$ then $\delta\theta_1^{\prime}=\pm \delta\theta_1,\delta\phi_1^{\prime}=\pm \delta\phi_1$.
Assume now that $\delta\theta_1\neq 0$ and $\delta\phi_1\neq 0$.
We obtain a contradiction since $a\cdot (\delta\theta_1,\delta\phi_1)=(\delta\theta_1,-\delta\phi_1)=\lambda(\delta\theta_1,\delta\phi_1)$ and $\lambda=\pm1$.
To end, we have to deal with polar tangent vectors:
$a\cdot (\delta x,\delta y,\delta x^{\prime},\delta y^{\prime})=(\delta x^{\prime},-\delta y^{\prime},\delta x,-\delta y)$ so $\delta x=\delta x^{\prime},\ \delta y=-\delta y^{\prime}$.
The action of $b$ leads to $\delta x=\delta x^{\prime}=\delta y=\delta y^{\prime}=0$.
It follows that the irreducible representations of dimension 1 are the ones listed in the proposition and only these ones.

Vectors $\alpha_{0,\phi}$ and $\alpha_{0,\theta}^{\prime}$ do not lie in $\G$ since  $D\Phi(x_e)\cdot \alpha_{0,\theta}^{\prime}=(0,0,-2N\sin(\theta_0))$ and $\alpha_{0,\phi}\in \so_\mu\cdot x_e$.
\ep
\begin{prop}
\label{rep2}
Let $q$ be an integer. 
Assume first $k_p=0$.
The following $\R$-spaces are irreducible representations of dimension $2$ of $G_{x_e}$ on $T_{x_e}\PP$:\\
\vspace{0.2cm}\\
$\star\; \<\alpha_{q,\theta},\beta_{q,\theta}\m ,\; \<\alpha_{q,\phi},\beta_{q,\phi}\m ,\; \<\alpha_{q,\theta}^{\prime},\beta_{q,\theta}^{\prime}\m ,\; \<\alpha_{q,\phi}^{\prime},\beta_{q,\phi}^{\prime}\m $, $1\leq q\leq N\! -\! 1$,  if  $N$ is odd and $G_{x_e}=D_{Nd}$ or $D_{Nh}$\\
\vspace{0.2cm}\\
$\star\; \<\alpha_{q,\theta},\beta_{q,\theta}\m ,\; \<\alpha_{q,\phi},\beta_{q,\phi}\m ,\; \<\alpha_{q,\theta}^{\prime},\beta_{q,\theta}^{\prime}\m ,\; \<\alpha_{q,\phi}^{\prime},\beta_{q,\phi}^{\prime}\m$, $1\leq q\leq N\! -\! 1$, $q\neq \frac{N}{2}$, $\<\alpha_{\frac{N}{2},\theta}^{\prime},\beta_{\frac{N}{2},\phi}\m,\; \<\alpha_{\frac{N}{2},\phi},\beta_{\frac{N}{2},\theta}^{\prime}\m$   if  $N$ is even and $G_{x_e}=D_{Nd}$\\
\vspace{0.2cm}\\
$\star\; \<\alpha_{q,\theta},\beta_{q,\theta}\m ,\; \<\alpha_{q,\phi},\beta_{q,\phi}\m,\; \<\alpha_{q,\theta}^{\prime},\beta_{q,\theta}^{\prime}\m,\; \<\alpha_{q,\phi}^{\prime},\beta_{q,\phi}^{\prime}\m$, $1\leq q\leq N\! -\! 1$, $q\neq\frac{N}{2}$, $\<\alpha_{\frac{N}{2},\theta},\alpha_{\frac{N}{2},\theta}^{\prime}\m,\; \<\alpha_{\frac{N}{2},\phi},\alpha_{\frac{N}{2},\phi}^{\prime}\m$  if  is $N$ even and $G_{x_e}=D_{Nh}$.\\
\vspace{0.2cm}\\
If $k_p=2$, then there are two supplementary irreducible representations of dimension $2$: $\<\delta x_p,\delta y_p\m$ and $\<\delta x_p^{\prime},\delta y_p^{\prime}\m$.\\
The representations $\<\alpha_{1,\theta},\beta_{1,\theta} \m$, $\<\alpha_{1,\phi}^{\prime},\beta_{1,\phi}^{\prime}\m$ and $\<\delta x_p^{\prime},\delta y_p^{\prime}\m$ do not lie in the symplectic slice, while the others do.
\end{prop}
\proof
Denote $V=T_{x_e}\PP$ and $V_{\C}=T_{x_e}\p\oplus iT_{x_e}\p$.
Irreducible representations of dimension 2 are $G_{x_e}$-invariant $\R$-planes of $V$.
We can identify $\R$-subspaces of dimension $\leq 2$ of $V$ with $\C$-lines of $V_{\C}$ by  $\< x,y \m_{\R} \mapsto \<z=x+iy\m_{\C}$ where $x,y\in V$.
It follows that $\< x_0,y_0 \m_{\R}$ is a $G_{x_e}$-invariant $\R$-plane of $V$ iff $x_0$ and $y_0$ are not colinear and $G_{x_e}$ acts like the transformation group $\< \{z\mapsto \lambda z,\lambda z\in \C  \}, z\mapsto \overline{z}\m$ on $\< z_0 \m_{\C}$ where  $z_0=x_0+iy_0\in V_{\C}$, that is $\< z_0 \m_{\C}\cup\< \overline{z_0}\m_{\C}$ is a $G_{x_e}$-invariant set of $V_{\C}$.

We have $G_{x_e}=D_{Nd}$ or $D_{Nh}$.
Isometries in question act on $V$ like:
$$
\left\lbrace\begin{array}{l}
s_x\cdot (\delta\theta_i,\delta\phi_i)=(\delta\theta_i,-\delta\phi_i)\\
s_z\cdot (\delta\theta_i,\delta\phi_i)=(-\delta\theta_i,\delta\phi_i)\\
r_{\alpha}\cdot (\delta\theta_i,\delta\phi_i)=(\delta\theta_i,\delta\phi_i)\\
s_x\cdot (\delta x,\delta y)=(\delta x,-\delta y)\\
s_z\cdot (\delta x,\delta y)=(\delta x,\delta y)\\
r_{\alpha}\cdot (\delta x,\delta y)=(\cos\alpha\;\delta x-\sin\alpha\;\delta y,\sin\alpha\;\delta x+\cos\alpha\;\delta y). 
\end{array}\right.
$$
Taking account of permutation symmetries, we obtain that the spaces listed in the proposition satisfy the above conditions.
These representations of dimension 2 are irreducible because they do not contain representations of dimension 1 which are listed in Prop. \ref{rep1}.
We then use Section \ref{sympslice} to examine wether these representations lie in the symplectic slice.
\ep

The following theorem lists the appropriate decompositions of the symplectic slice.
\begin{theo}
\label{mainth}
Let $N \geq 3$ and $x_e$ be of type $D_{Nd}(R,R',k_pp)$ or $D_{Nh}(2R,k_pp)$ such that $\mu=\Phi(x_e)$ is non-zero.
Let $l$ be an integer and denote 
$$B^0=\<\alpha_{0,\theta},\alpha_{0,\phi}^{\prime},\beta_{1,\phi},\alpha_{1,\theta}^{\prime},\alpha_{1,\phi},\beta_{1,\theta}^{\prime}\m$$
$$B_p^0=\<\alpha_{0,\theta},\alpha_{0,\phi}^{\prime},\beta_{1,\phi},\alpha_{1,\theta}^{\prime},\delta x_p,\alpha_{1,\phi},\beta_{1,\theta}^{\prime},\delta y_p\m$$
$$B^1=\<\sin\theta_0\;\alpha_{1,\theta}+\cos\theta_0\;\beta_{1,\phi}^{\prime},\sin\theta_0\;\beta_{1,\theta}+\cos\theta_0\;\alpha_{1,\phi}^{\prime}\m$$
$$B_p^1=B^1\oplus\<\alpha_{1,\theta}-\frac{N}{2}\cos\theta_0\;\delta x_p^{\prime},\beta_{1,\theta}-\frac{N}{2}\cos\theta_0\;\delta y_p^{\prime}\m$$
$$B_l=\<\{\alpha_{q,\theta},\beta_{q,\phi}^{\prime},\beta_{q,\theta},\alpha_{q,\phi}^{\prime},\alpha_{q,\theta}^{\prime},\beta_{q,\phi},\beta_{q,\theta}^{\prime},\alpha_{q,\phi}\mid 2\leq q\leq l\}\m.$$  
Let $k_p=0$. 
With the following decompositions of the symplectic slice $\G$, the Hessian  $d^2H_{\xi}|_{\G}(x_e)$ diagonalizes in $2\times 2$ blocks  with four  $1\times 1$ blocks:\\
\vspace{0.2cm}\\
$\star\; B^0\oplus B^1\oplus B_{\frac{N-1}{2}}\; $ for $x_e=D_{Nd}(R,R')$ or $D_{Nh}(2R)$, and $N$ odd\\
\vspace{0.2cm}\\
$\star\; B^0\oplus B^1\oplus B_{\frac{N}{2}\!-\!1}\oplus \<\alpha_{\frac{N}{2},\theta}^{\prime},\beta_{\frac{N}{2},\phi},\beta_{\frac{N}{2},\theta}^{\prime},\alpha_{\frac{N}{2},\phi}\m$ for $x_e=D_{Nd}(R,R')$ and $N$ even\\
\vspace{0.2cm}\\
$\star\; B^0\oplus B^1\oplus B_{\frac{N}{2}\!-\!1}\oplus \<\alpha_{\frac{N}{2},\theta},\alpha_{\frac{N}{2},\theta}^{\prime},\alpha_{\frac{N}{2},\phi},\alpha_{\frac{N}{2},\phi}^{\prime}\m$ for $x_e=D_{Nh}(2R)$ and $N$ even.\\
\vspace{0.2cm}\\
\noindent Let $k_p=2$. 
With the following decompositions of $\G$, the Hessian  $d^2H_{\xi}|_{\G}(x_e)$ diagonalizes in  $2\times 2$ blocks  with two  $3\times 3$ blocks  and two  $1\times 1$ blocks:\\
\vspace{0.2cm}\\
$\star\; B_p^0\oplus B_p^1\oplus B_{\frac{N-1}{2}}\; $ for $x_e=D_{Nd}(R,R',2p)$ or $D_{Nh}(2R,2p)$, and $N$ odd\\
\vspace{0.2cm}\\
$\star\; B_p^0\oplus B_p^1\oplus B_{\frac{N}{2}\!-\!1}\oplus \<\alpha_{\frac{N}{2},\theta}^{\prime},\beta_{\frac{N}{2},\phi},\beta_{\frac{N}{2},\theta}^{\prime},\alpha_{\frac{N}{2},\phi}\m$ for $x_e=D_{Nd}(R,R',2p)$ and $N$ even\\
\vspace{0.2cm}\\
$\star\; B_p^0\oplus B_p^1\oplus B_{\frac{N}{2}\!-\!1}\oplus \<\alpha_{\frac{N}{2},\theta},\alpha_{\frac{N}{2},\theta}^{\prime},\alpha_{\frac{N}{2},\phi},\alpha_{\frac{N}{2},\phi}^{\prime}\m$ for $x_e=D_{Nh}(2R,2p)$ and $N$ even.\\
\vspace{0.2cm}\\
\noindent With these decompositions, the matrix $L=\Omega_{\G}^{\flat -1}d^2H_{\xi}|_{\G}(x_e)$ of the linearized system in the symplectic slice (where $\Omega_{\G}^{\flat}$ is the matrix of $\omega|_{\G}$) diagonalizes in $4\times 4$ blocks  with two  $2\times 2$ blocks if $k_p=0$, and in  $4\times 4$ blocks  with one  $6\times 6$ block  and one  $2\times 2$ block  if $k_p=2$.
\end{theo}

The case $N=2$ will be treated case by case later.\\
\\
\proof
There are no difficulties to check that vectors listed form a basis of the symplectic slice using Prop. \ref{diffmom}, \ref{tgtspace}, \ref{rep1},  \ref{rep2} and keeping in mind that $\dim \G=4N-4$ (resp. $4N$) if $k_p=0$ (resp. $k_p=2$).

Let $x_e$ be a $D_{Nd}(R,R',k_pp)$ arrangement.
We have $G_{x_e}=C_{2N}\rtimes \Z_2[s_x]$ where $s_x$ is the reflection $x\mapsto -x$.
Note that $C_{2N}$ is symplectic while $s_x$ is anti-symplectic.\\

Assume $k_p=0$.
First, one examines whether some irreducible representations of $C_{2N}$ of the symplectic slice $\G$ are isomorphic.
The two irreducible representations of $C_{2N}$ of dimension 1 are isomorphic.
For irreducible representations of $C_{2N}$ of dimension 2, we obtain couples of isomorphic representations:
subspaces $\<\alpha_{q,\theta},\beta_{q,\phi}^{\prime},\beta_{q,\theta},\alpha_{q,\phi}^{\prime},\alpha_{q,\theta}^{\prime},\beta_{q,\phi},\beta_{q,\theta}^{\prime},\alpha_{q,\phi}\m$, $q\geq 2$ are split in $\<\alpha_{q,\theta},\beta_{q,\phi}^{\prime},\beta_{q,\theta},\alpha_{q,\phi}^{\prime}\m$ and $\<\alpha_{q,\theta}^{\prime},\beta_{q,\phi},\beta_{q,\theta}^{\prime},\alpha_{q,\phi}\m$ which correspond respectively to rotations of angle ${q\pi}/{N}$ and $-{q\pi}/{N}$.
The representation $B^1$ is not isomorphic to another one, while $\<\beta_{1,\phi},\alpha_{1,\theta}^{\prime}\m$ and $\<\alpha_{1,\phi},\beta_{1,\theta}^{\prime}\m$ are isomorphic.
Hence, we get the isotypic decomposition of $\G$, say $\G=\oplus_j V_j$ with $V_0=\<\alpha_{0,\theta},\alpha_{0,\phi}^{\prime}\m$, $V_1=B^1$ and $\dim V_j=4$ ($j\geq 2$).
Since $H$ is $G_{x_e}$-invariant, it follows that $d^2H_{\xi}(x_e)$ is $G_{x_e}$-invariant and the linear map $\widetilde{d^2H_{\xi}|_{\G}(x_e)}:\oplus_j V_j \to \oplus_j V_j^{\ast}$ is $G_{x_e}$-equivariant.
Let $i_{V_j}: V_j \to \G$ and $p_{V_j^{\ast}}: \G^{\ast} \to V_j^{\ast}$ be respectively the canonical inclusion and projection.
By Shur's Lemma, $$p_{V_j^{\ast}}\circ\widetilde{d^2H_{\xi}|_{\G}(x_e)}\circ i_{V_i}=0\ \mbox{for all}\ i\neq j$$ 
(representations $V_j$ and $V_j^{\ast}$ are isomorphic since $G_{x_e}$ is compact).
Since the dynamical system is $C_{2N}$-equivariant, this result holds also for the linear map $L=\Omega_{\G}^{\flat -1}d^2H_{\xi}|_{\G}(x_e)$: $L$ diagonalizes in $4\times 4$ blocks  with two $2\times 2$ blocks.

The linear map $p_{V_j^{\ast}}\circ\widetilde{d^2H_{\xi}|_{\G}(x_e)}\circ i_{V_j}$ is $\Z_2[s_x]$-equivariant and subspaces $V_j$ split in two subspaces of dimension $\frac{1}{2}\dim V_j$ (the $\Z_2[s_x]$-isotypic decomposition of $V_j$).
That is, $V_j=U_j^0\oplus U_j^1$ and $$p_{U_j^{k\ast}}\circ\widetilde{d^2H_{\xi}|_{\G}(x_e)}\circ i_{U_j^l}=0\ \mbox{for all}\ k\neq l.$$
The Hessian diagonalizes therefore in $2\times 2$ blocks  with four  $1\times 1$ blocks.\\

Assume $k_p=2$.
With respect to the case $k_p=0$, the changes are:\\
$\bullet$ $B^1$ is isomorphic to the representation $\<\alpha_{1,\theta}-\frac{N}{2}\cos\theta_0\;\delta x_p^{\prime},\beta_{1,\theta}-\frac{N}{2}\cos\theta_0\;\delta y_p^{\prime}\m$, this provides a $4\times 4$ block in $L$ and two $2\times 2$ blocks in the Hessian.\\
$\bullet$ we have a triplet of isomorphic representations 
$$\<\beta_{1,\phi},\alpha_{1,\theta}^{\prime}\m,\ \<\alpha_{1,\phi},\beta_{1,\theta}^{\prime}\m\ \mbox{and}\ \<\delta x_p,\delta y_p\m$$ 
which provides a $6\times 6$ block in $L$ and two $3\times 3$ blocks in the Hessian.\\

For $D_{Nh}(2R,k_pp)$ arrangements, $G_{x_e}=C_{N}\times\Z_2[(s_z,\tau)]\rtimes \Z_2[s_x]$.
The proof is similar: $C_{N}\times\Z_2[(s_z,\tau)]$ is symplectic while  $\Z_2[s_x]$ is anti-symplectic.        \ep
    
\rmk 
The appropriate basis for the block diagonalization does not depend on the expression for the Hamiltonian.
Let $M$ be the product of $2N$ spheres $S^2$, $x\in M$ with symmetry group $D_{Nd}$ or $D_{Nh}$, and $\f: M \to \R$ a $G_x$-invariant map.
The Hessian of $\f$ at $x$ block diagonalizes like in the case $k_p=0$ of the previous Theorem.

\subsubsection{Expression of $d^2H_{\xi}(x_e)$}
\label{hessian}

First, consider the case $k_p=0$ (no polar vortices).
Recall that the Hamiltonian is 
$$H=\sum_{i<j} \lambda_i  \lambda_j \ln (\, 2(1-\cos\theta_i \cos\theta_j -\sin\theta_i \sin\theta_j \cos(\phi_i-\phi_j))\,)$$
and $H_{\xi}=H+\xi\sum_{i=1}^{2N}\lambda_i\cos\theta_i$.\\
Relative equilibria $D_{Nd}(R,R')$ and $D_{Nh}(2R)$ are characterized by coordinates:
$$  
 \left\lbrace
            \begin{array}{l}
                \theta_j=\theta_0\ \mbox{for}\ j=1,\dots,N\\
                \theta_j^{\prime}=\pi-\theta_0\ \mbox{for}\ j=1,\dots,N
            \end{array}
       \right.
$$
$$
\mbox{and}\     \left\lbrace
            \begin{array}{l}
                \phi_j=\frac{2\pi(j-1)}{N}\ \mbox{for}\ j=1,\dots,N\\
                \phi_j^{\prime}=\left\lbrace \begin{array}{l}
                     \frac{2\pi(j-1)}{N}\ \mbox{if}\ x_e=D_{Nh}(2R)\\
                     \frac{2\pi(j-1)}{N}+\frac{\pi}{N}\ \mbox{if}\ x_e=D_{Nd}(R,R')                                                      \end{array} \right.
                                                   \mbox{for}\ j=1,\dots,N.
            \end{array}
       \right.
$$
Note that we assume $\theta_0\in ]0,\frac{\pi}{2}]$ without loss of generality since the relative equilibria in question have a ``North-South'' symmetry.
From now on, we denote $u=\cos\theta_0$, $\varphi_j=\frac{2\pi(j-1)}{N}$ and $\varphi_j^{\prime}=\varphi_j$ (resp. $\varphi_j+\frac{\pi}{N}$) if $x_e$ is a $D_{Nh}$ (resp. $D_{Nd}$) relative equilibrium, in order to lighten formula.

Some tedious computations show that second derivarives of $H_{\xi}$ at $x_e$ are:
$$ 
\left\lbrace\begin{array}{l}
\frac{\partial^2 H_{\xi}}{\partial \theta_i^2}=\frac{\partial^2  H_{\xi}}{\partial \theta_i^{\prime 2}}=\frac{1}{1-u^2}\sum_{j=2}^{N}\frac{\cos\varphi_j}{1-\cos\varphi_j} \\
\ \ \ \ \ \ \ \ \ \ \ \ \ \ \ \ \ -\sum_{j=1}^{N}\frac{-2u^2+(1-u^2)\cos\varphi_j^{\prime}-(1-u^2)\cos^2\varphi_j^{\prime}}{[1+u^2-\cos\varphi_j^{\prime}\; (1-u^2)]^2}  -\xi u\\
\frac{\partial^2 H_{\xi}}{\partial\theta_i\partial\theta_j}=\frac{\partial^2 H_{\xi}}{\partial\theta_i^{\prime}\partial\theta_j^{\prime}}=\frac{-1}{(1-u^2)(1-\cos\varphi_{j\!-\!i\!+\!1})},\ i\neq j\\
\frac{\partial^2 H_{\xi}}{\partial\theta_i\partial\theta_j^{\prime}}=\frac{1-u^2-\cos\varphi_{j\!-\!i\!+\!1}^{\prime}\; (1+u^2)}{[1+u^2-\cos\varphi_{j\!-\!i\!+\!1}^{\prime}\; (1-u^2)]^2}\\
\frac{\partial^2 H_{\xi}}{\partial\theta_i\partial\phi_j}=\frac{\partial^2 H_{\xi}}{\partial\theta_i^{\prime}\partial\phi_j^{\prime}}=0\\
\frac{\partial^2 H_{\xi}}{\partial\theta_i\partial\phi_j^{\prime}}=\frac{\partial^2 H_{\xi}}{\partial\phi_i\partial\theta_j^{\prime}}=-2u\sqrt{1-u^2}\frac{\sin\varphi_{j\!-\!i\!+\!1}^{\prime}}{[1+u^2-\cos\varphi_{j\!-\!i\!+\!1}^{\prime}\; (1-u^2)]^2}\\
\frac{\partial^2 H_{\xi}}{\partial\phi_i^2}=\frac{\partial^2 H_{\xi}}{\partial\phi_i^{\prime 2}}=\sum_{j=2}^{N}\frac{-1}{1-\cos\varphi_j} +(1-u^2)\sum_{j=1}^{N}\frac{1-u^2-\cos\varphi_j^{\prime}\; (1+u^2)}{[1+u^2-\cos\varphi_j^{\prime}\; (1-u^2)]^2}\\
\frac{\partial^2 H_{\xi}}{\partial\phi_i\partial\phi_j}=\frac{\partial^2 H_{\xi}}{\partial\phi_i^{\prime}\partial\phi_j^{\prime}}=\frac{1}{1-\cos\varphi_{j\!-\!i\!+\!1}},\ i\neq j\\
\frac{\partial^2 H_{\xi}}{\partial\phi_i\partial\phi_{j}^{\prime}}=-(1-u^2)\frac{1-u^2-\cos\varphi_{j\!-\!i\!+\!1}^{\prime}\; (1+u^2)}{[1+u^2-\cos\varphi_{j\!-\!i\!+\!1}^{\prime}\; (1-u^2)]^2}
\end{array}
\right.    
$$
where the sum in the subscripts is modulo $N$.

Consider now the second case: $k_p=2$ (with polar vortices).
If we permit the co-latitude of the relative equilibria $\theta_0$ to lie in $]0,\pi[$ instead of $]0,\frac{\pi}{2}]$, we can choose $\lambda_n=+1$ and $\lambda_s=-1$ without loss of generality.
The Hamiltonian becomes:
$$
H=\sum_{i<j}^{2N} \lambda_i  \lambda_j \ln (\, 2(1-\cos\theta_i \cos\theta_j -\sin\theta_i \sin\theta_j \cos(\phi_i-\phi_j))\,)+H_{r,n}+H_{r,s}+H_{n,s}
$$
where $$H_{r,n}=\sum_{i=1}^{2N} \lambda_i   \ln (\, 2(1-\sin\theta_i \cos\phi_i\; x_n -\sin\theta_i \sin\phi_i\; y_n-\sqrt{1-x_n^2-y_n^2}\cos\theta_i)\,)$$
$$H_{r,s}= -\sum_{i=1}^{2N} \lambda_i   \ln (\, 2(1-\sin\theta_i \cos\phi_i\; x_s -\sin\theta_i \sin\phi_i\; y_s+\sqrt{1-x_s^2-y_s^2}\cos\theta_i)\,)$$
$$H_{n,s}= -\ln (\, 2(1-x_n x_s - y_n y_s+\sqrt{1-x_n^2-y_n^2}\sqrt{1-x_s^2-y_s^2})\,)$$
and augmented Hamiltonian:
$$H_{\xi}=H+\xi\left(\sum_{i=1}^{2N}\lambda_i\cos\theta_i +\sqrt{1-x_n^2-y_n^2}-\sqrt{1-x_s^2-y_s^2}\right).$$
Note that the angular velocity is different than in the case $k_p=0$ (see Prop. \ref{angvelexpr}).
Second derivatives of $H_{\xi}$ at $x_e$ are those given above together with:
$$
\left\lbrace\begin{array}{ll} 
\frac{\partial^2 H_{\xi}}{\partial\theta_i\partial x_n}=\frac{\cos\varphi_i}{1-u}
& \frac{\partial^2 H_{\xi}}{\partial\theta_i^{\prime}\partial x_n}=-\frac{\cos\varphi_i^{\prime}}{1+u}\\
\frac{\partial^2 H_{\xi}}{\partial\theta_i\partial x_s}=\frac{\cos\varphi_i}{1+u}
& \frac{\partial^2 H_{\xi}}{\partial\theta_i^{\prime}\partial x_s}=-\frac{\cos\varphi_i^{\prime}}{1-u}\\
\frac{\partial^2 H_{\xi}}{\partial\theta_i\partial y_n}=\frac{\sin\varphi_i}{1-u}
& \frac{\partial^2 H_{\xi}}{\partial\theta_i^{\prime}\partial y_n}=-\frac{\sin\varphi_i^{\prime}}{1+u}\\
\frac{\partial^2 H_{\xi}}{\partial\theta_i\partial y_s}=\frac{\sin\varphi_i}{1+u}
& \frac{\partial^2 H_{\xi}}{\partial\theta_i^{\prime}\partial y_s}=-\frac{\sin\varphi_i^{\prime}}{1-u}\\
\frac{\partial^2 H_{\xi}}{\partial\phi_i\partial x_n}=\frac{\sqrt{1-u^2}}{1-u}\sin\varphi_i
& \ \frac{\partial^2 H_{\xi}}{\partial\phi_i^{\prime}\partial x_n}=-\frac{\sqrt{1-u^2}}{1+u}\sin\varphi_i^{\prime}\\
\frac{\partial^2 H_{\xi}}{\partial\phi_i\partial x_s}=-\frac{\sqrt{1-u^2}}{1+u}\sin\varphi_i
& \frac{\partial^2 H_{\xi}}{\partial\phi_i^{\prime}\partial x_s}=\frac{\sqrt{1-u^2}}{1-u}\sin\varphi_i^{\prime}\\
\frac{\partial^2 H_{\xi}}{\partial\phi_i\partial y_n}=-\frac{\sqrt{1-u^2}}{1-u}\cos\varphi_i
& \frac{\partial^2 H_{\xi}}{\partial\phi_i^{\prime}\partial y_n}=\frac{\sqrt{1-u^2}}{1+u}\cos\varphi_i^{\prime}\\
\frac{\partial^2 H_{\xi}}{\partial\phi_i\partial y_s}=\frac{\sqrt{1-u^2}}{1+u}\cos\varphi_i
& \frac{\partial^2 H_{\xi}}{\partial\phi_i^{\prime}\partial y_s}=-\frac{\sqrt{1-u^2}}{1-u}\cos\varphi_i^{\prime}\\
\frac{\partial^2 H_{\xi}}{\partial x_n^2}=\frac{\partial^2 H_{\xi}}{\partial x_s^2}=\frac{\partial^2 H_{\xi}}{\partial y_n^2}=\frac{\partial^2 H_{\xi}}{\partial y_s^2}=1/2-\xi
& \frac{\partial^2 H_{\xi}}{\partial x_n\partial x_s}=\frac{\partial^2 H_{\xi}}{\partial y_n\partial y_s}=1/2\\
\frac{\partial^2 H_{\xi}}{\partial x_n\partial y_n}=\frac{\partial^2 H_{\xi}}{\partial x_n\partial y_s}=\frac{\partial^2 H_{\xi}}{\partial x_s\partial y_s}=\frac{\partial^2 H_{\xi}}{\partial x_s\partial y_n}=0.
\end{array}\right.
$$

\subsubsection{Stability results}

One is able now to obtain an expression of $d^2H_{\xi}|_{\G}(x_e)$ and $L=\Omega_{\G}^{\flat -1}d^2H_{\xi}|_{\G}(x_e)$ formed of diagonal blocks.
We shall see in the proofs of this section that this block diagonalization permits us to derive some formulae for the eigenvalues of $d^2H_{\xi}|_{\G}(x_e)$ and $L$, and thus to conclude about both Lyapunov and linear stability.
We first study the linearization since linear instability is easier to prove than Lyapunov stability.
Indeed, to prove linear instability, it is sufficient that one block is linearly unstable, while to prove Lyapunov stability, one needs that the blocks are all positive definite or all negative definite.
In the statements of this section, it is understood that when a \rum $x$ is claimed to be \emph{linearly stable}, it means also that the Hessian at $x$ is not definite (thus we can not conclude on Lyapunov stability applying the energy momentum method).
Recall that a linearly \emph{stable} \rum can be Lyapunov \emph{unstable} (see \cite{MR94} for an example).
The following theorem gives results for $D_{Nd}(R,R')$ relative equilibria.
\begin{theo}
\label{DNd}
\ 
 
\noindent $\bullet$ $D_{2d}(R,R')$ relative equilibria are Lyapunov stable for $\theta_0\in[1.14,\pi/2]$ and linearly unstable on the complement in $]0,\pi/2]$.\\
$\bullet$ $D_{3d}(R,R')$ relative equilibria are linearly stable on $[1.302,1.315]$ and linearly unstable on the complement in $]0,\pi/2]$.
On the upper bound occurs an Hamilton-Hopf bifurcation.\\
$\bullet$ For $N\geq 4$, $D_{Nd}(R,R')$ relative equilibria are linearly unstable for all co-latitude $\theta_0$.\\
$\bullet$ Equatorial fixed equilibria $D_{2Nh}(R_e)$ are linearly unstable for all $N\geq 3$, while $D_{4h}(R_e)$ is Lyapunov stable.
\end{theo}
\proof
Let $x_e$ be a $D_{Nd}(R,R')$ relative equilibria and $N\geq 3$.
By Theorem \ref{mainth}, the linearized system $L=\Omega^{\flat -1}d^2H_{\xi}|_{\G}(x_e)$ diagonalizes in $4\times 4$ blocks and two $2\times 2$ blocks.
All diagonal $4\times 4$ blocks of $d^2H_{\xi}|_{\G}(x_e)$ and $J$ are respectively like (up to a factor):
$$
\left( \begin{array}{cccc} 
a & c & 0 & 0\\
c & b & 0 & 0\\
0 & 0 & a & -c\\
0 & 0 & -c & b
\end{array}\right)
\ ,\ \left( \begin{array}{cccc} 
0 & 0 & 0 & -1\\
0 & 0 & 1 & 0\\
0 & -1 & 0 & 0\\
1 & 0 & 0 & 0
\end{array}\right) .
$$
Hence, diagonal $4\times 4$ blocks of $L$ are like (up to a factor):  
$$
\left( \begin{array}{cccc} 
0 & 0 & c & -b\\
0 & 0 & a & -c\\
-c & -b & 0 & 0\\
a & c & 0 & 0
\end{array}\right)
$$
whose eigenvalues are $\pm i(c+\sqrt{ab}),\pm i(c-\sqrt{ab})$.
These eigenvalues have a non-zero real part if and only if $ab<0$.
Note that the matrix above is equal to $cD+aN_1+bN_2$ where $D$ is semi-simple and $N_1,N_2$ are nilpotent.
When $ab$ passes through zero, linear stability is lost through a Hamiltonian-Hopf bifurcation if $c\neq 0$.
Note also that the corresponding block of the Hessian is definite if and only if $ab>c^2$. 

The two $2\times 2$ blocks of $L$ correspond to the subspaces $\<\alpha_{0,\theta},\alpha_{0,\phi}^{\prime}\m$ and $B^1$ (see Theorem \ref{mainth}), they are respectively:
$$
\left( \begin{array}{cc} 
0 & -s\\
r & 0\\
\end{array}\right)
\ ,\ \left( \begin{array}{cc} 
0 & -w\\
w & 0\\
\end{array}\right)
$$
$$ 
\mbox{where}\ \left\lbrace\begin{array}{lll}
r=d^2H_{\xi}(x_e)\cdot (\alpha_{0,\tet},\alpha_{0,\tet})\\
s=d^2H_{\xi}(x_e)\cdot (\alpha_{0,\phi}^{\prime},\alpha_{0,\phi}^{\prime})\\
w=d^2H_{\xi}(x_e)\cdot (\sin\theta_0\;\alpha_{1,\theta}+\cos\theta_0\;\beta_{1,\phi}^{\prime},\sin\theta_0\;\alpha_{1,\theta}+\cos\theta_0\;\beta_{1,\phi}^{\prime}).
\end{array}\right. 
$$ 
The eigenvalues of these blocks are  respectively $\pm i\sqrt{rs}$ and $\pm iw$.
It follows that the subspace $B^1$ is linearly stable and that bifurcation occurs when $rs$ passes through zero.\\

We first deal with the $2\times 2$ block of $L$ corresponding to $\<\alpha_{0,\theta},\alpha_{0,\phi}^{\prime}\m$ .
Some lengthy computations show that the expressions of $r$ and $s$ are:
$$ 
\begin{array}{ll}
r=2N \left[ -(N-1)\frac{1+u^2}{1-u^2}+(1-u^4)\Sigma_1-2(1+u^4)\Sigma_2+(1-u^4)\Sigma_3 \right]\\
s=4N(1-u^2) \left[ (1-u^2)\Sigma_1-(1+u^2)\Sigma_2  \right]
\end{array}
$$ 
where $\Sigma_m=\sum_{j=1}^N\frac{\cos^{m-1} \varphi_j^{\prime}}{[1+u^2-(1-u^2)\cos \varphi_j^{\prime}]^2}$.
When $N$ goes to infinity, $r/N^2$ goes to $0$ and $s/N^2$ to $-\frac{2u}{1-u^2}$.
One checks that $r>0$ for all $\theta_0\in [0,{\pi}/{2}]$.
One can prove that $s$ increases strictly with respect to $\theta_0$.
Since limits of $s$ at $0$ and ${\pi}/{2}$ are respectively negative and positive, relative equilibria are linearly unstable for $\theta_0\in ]0,t_N[$ where $t_N$ is defined by $s(t_N)=0$.
When $N$ goes to infinity, $t_N$ goes to ${\pi}/{2}$.

Let $N\geq4$.
Consider the blocks of $L$ corresponding to the subspaces $\<\alpha_{q,\theta}^{\prime},\beta_{q,\phi},\beta_{q,\theta}^{\prime},\alpha_{q,\phi}\m$. 
Following the notations of the beginning of the proof, we have:
$$ 
\left\lbrace\begin{array}{ll}
a=d^2H_{\xi}(x_e)\cdot (\alpha_{q,\tet}^{\prime},\alpha_{q,\tet}^{\prime})\\
b=d^2H_{\xi}(x_e)\cdot (\beta_{q,\phi},\beta_{q,\phi})\\
c=d^2H_{\xi}(x_e)\cdot (\alpha_{q,\tet}^{\prime},\beta_{q,\phi}).
\end{array}\right.
$$ 
Some lengthy computations show that:
$$
\begin{array}{lll}
a & = & N [ -(N-1)\frac{u^2}{1-u^2}+\frac{1}{1-u^2}\sum_{j=2}^N\frac{\cos \varphi_j-\cos q\varphi_j}{1-\cos \varphi_j}\\[10pt]
  & - & \sum_{j=1}^N  \frac{u^4-u^2+(1-u^2)\cos q\varphi_j^{\prime}+(2u^4-u^2+1-(1+u^2)\cos q\varphi_j^{\prime})\cos \varphi_j^{\prime}-(1-u^4)\cos^2 \varphi_j^{\prime}}{[1+u^2-(1-u^2)\cos \varphi_j^{\prime}]^2} ]\\[15pt]
b & = & N [ -\sum_{j=2}^N\frac{1-\cos q\varphi_j}{1-\cos \varphi_j}  +(1-u^2)\sum_{j=1}^N  (1-\cos q\varphi_j^{\prime})\frac{1-u^2-(1+u^2)\cos \varphi_j^{\prime}}{[1+u^2-(1-u^2)\cos \varphi_j^{\prime}]^2} ].
\end{array}
$$
One checks that if $q$ is the integer part of $N/2$, then $ab<0$ for $\theta_0\in [t_N,{\pi}/{2}]$ and for all $N\geq4$; $L$ has therefore eigenvalues with non-zero real part, and the relative equilibria are linearly unstable. 

For $N=3$, the matrix $L$ is formed of the two $2\times 2$ blocks with one $4\times 4$ block, the first $2\times 2$ block is linearly unstable iff $\theta_0\in [0,1.302]$ and the $4\times 4$ block is linearly unstable iff $\theta_0\in [1.315,\pi/2]$ (the $a$ and $b$ involved are given by the above formulae substituting $q=1$). 
Since the second block (subspace $B^1$) is linearly stable, $D_{3d}(R,R')$ is linearly stable on $[1.302,1.315]$. 
The Hessian is not definite on this interval, thus we can not know wether the \ria are Lyapunov stable on this interval.

The $D_{2Nh}(R_e)$ fixed equilibria correspond to $D_{Nd}(R,R')$ relative equilibria with $\theta_0={\pi}/{2}$.
In this case $\mu=0$ and $\so_{\mu}\cdot x_e$ is generated by $\alpha_{0,\phi}$, $\alpha_{1,\theta}$, and $\beta_{1,\theta}$ (Prop. \ref{tgtspace}). 
The symplectic slice is therefore $\G=\Ker D\Phi(x_e)\cap \<\alpha_{0,\phi},\alpha_{1,\theta},\beta_{1,\theta}\m^{\bot}$.
Hence eigenvalues corresponding to the subspace $B^1$ can be ignored when $\theta_0={\pi}/{2}$.
Consider the $4\times 4$ block corresponding to the subspace $\<\alpha_{1,\theta}^{\prime},\beta_{1,\phi},\beta_{1,\theta}^{\prime},\alpha_{1,\phi}\m$: 
for $\theta_0={\pi}/{2}$ and $q=1$, one has $a=-2N\sum_{j=1}^N\frac{\cos \varphi_j^{\prime}}{1-\cos \varphi_j^{\prime}}$ and $b=N$.   
Since $a<0$ for $N\geq 3$, the $D_{2Nh}(R_e)$ fixed equilibria are linearly unstable for $N\geq 3$.\\

Consider now the case $N=2$, that is the $D_{2d}(R,R')$ \riab. 
For $\mu\neq 0$, a basis of the symplectic slice is $$(\alpha_{0,\theta},\ \alpha_{0,\phi}^{\prime},\ \sin\theta_0\;\alpha_{1,\theta}-\cos\theta_0\;\beta_{1,\phi}^{\prime},\ \sin\theta_0\;\beta_{1,\theta}^{\prime}-\cos\theta_0\; \alpha_{1,\phi}).$$
The Hessian is a diagonal matrix in this basis.
A study of the Hessian and the linearized system gives desired results.
\ep

Then, consider the $D_{Nd}(R,R',2p)$ relative equilibria.
We recall that since we assume $\lambda_n=+1$, we allow $\theta_0$ to lie in $]0,\pi[$.
\begin{theo}
\label{DNdp}
\ 

\noindent $\bullet$ $D_{2d}(R,R',2p)$ relative equilibria are linearly stable on $[2.21,2.31]$ and linearly unstable on the complement in $]0,\pi[$.\\
$\bullet$ $D_{3d}(R,R',2p)$ relative equilibria are linearly stable on $[1.8,2.05]\cup[2.25,\pi[$ and linearly unstable on the complement in $]0,\pi[$.\\ 
$\bullet$ For $N\geq 4$, $D_{Nd}(R,R',2p)$ relative equilibria are linearly stable on a small interval $I_N$ and linearly unstable on the complement of this interval in $]0,\pi[$. 
For $N=4,\dots,7$, these intervals are: $I_4=[1.75,1.79]$, $I_5=[1.73,1.76]$, $I_6=[1.71,1.72]$, $I_7=[1.69,1.70]$; while the length of $I_N$ is less than $10^{-2}$ for $N\geq 8$ .
Hamiltonian-Hopf bifurcations occur at the bounds of $I_N$.
\end{theo}
\proof
The proof is similar to that for $D_{Nd}(R,R')$.
The components of the Hessian are equal to the components of the case $k_p=0$ except the components $d^2H_{\xi}(x_e)\cdot (\delta \theta_i,\delta \theta_i)$ and of course components involving tangent vectors at the poles (see Section \ref{hessian}).
Hence, expressions of $s$ and $b$ do not change, while a term $-\frac{4Nu}{1-u^2}\lambda_n$ is added to $a$ and a term $-\frac{8Nu}{1-u^2}\lambda_n$ is added to $r$.
The subspace $B_p^1$ is no longer linearly stable. 
The expressions of $a$ and $b$ when $q=1$ do not correspond to any eigenvalue formula since $\alpha_{1,\theta}^{\prime}$ and $\beta_{1,\phi}$ lie in the subspace $\<\beta_{1,\phi},\alpha_{1,\theta}^{\prime},\delta x_p,\alpha_{1,\phi},\beta_{1,\theta}^{\prime},\delta y_p\m$ which correspond to the $6\times 6$ block.
Thus the expressions of $a$ and $b$ are useful only for $N\geq 4$.

For $N\geq 4$, the linear stability is governed by the block corresponding to the subspace $\<\alpha_{q,\theta}^{\prime},\beta_{q,\phi},\beta_{q,\theta}^{\prime},\alpha_{q,\phi}\m$ where $q$ is the integer part of ${N}/{2}$ (as in the case $k_p=0$). 
Hamiltonian-Hopf bifurcations occur when linear stability is lost (see the previous proof).

The stability results hold also when $\mu=0$:
to obtain the symplectice slice when $\mu$ vanishes, one has to remove the subspace $$Q=\< \beta_{1,\theta}+\cos\theta_0\sin\theta_0\;\alpha_{1,\phi}^{\prime}+\delta y_p^{\prime},\ \alpha_{1,\theta}-\cos\theta_0\sin\theta_0\;\beta_{1,\phi}^{\prime}+\delta x_p^{\prime}\m$$ from the symplectice slice when $\mu\neq 0$ (Prop. \ref{tgtspace}).
For $N\geq 4$, the subspace governing stability has an empty intersection with $Q$, so the stability of the zero-momentum configuration is the stability of this configuration found with the generic symplectic slice (where $Q$ is non removed).
For $N=2,3$, one checks by hand that stability results hold also when $\mu=0$.
\ep

\rmk
Adding two polar vortices such that the ``minus'' polar vortex and the ``plus'' ring are in the same hemisphere, tends to stabilize $D_{Nd}(R,R')$ relative equilibria (above all for $N=3$):
in the sense that the stronger are the polar vortices, the larger is the domain of stability.\\

The following theorem lists stability results for $D_{Nh}(2R)$ relative equilibria.
\begin{theo}
\label{DNh}
\
 
\noindent $\bullet$ For $N \geq 7$, $D_{Nh}(2R)$ relative equilibria are linearly unstable for all co-latitude $\theta_0$.\\
$\bullet$ The following relative equilibria are Lyapunov stable in the given range ($\theta_0\in ]0,\pi/2[$) and linearly unstable in the complement:\\
$D_{2h}(2R)$  $]0,0.66]$\\
$D_{4h}(2R)$  $]0,0.73]$\\
$D_{6h}(2R)$  $]0,0.45]$.\\
A pitchfork bifurcation occurs at the upper bound of these intervals.\\
$\bullet$ The following relative equilibria are Lyapunov stable on the first range, linearly stable on the second range and linearly unstable on the third one:\\
$D_{3h}(2R)$  $]0,0.77]$, $]0.77,0.78]$, $]0.78,\pi/2]$\\
$D_{5h}(2R)$  $]0,0.67]$, $]0.67,0.68]$, $]0.68,\pi/2]$.\\
Linear stability is lost through a Hamiltonian-Hopf bifurcation.
\end{theo}
\proof
Assume first that $N$ is odd.
The basis for the symplectic slice is the same than for $D_{Nd}(R,R')$ (Theorem \ref{mainth}).
The eigenvalues of the linearized system $L$ on the subspace $\<\alpha_{q,\theta}^{\prime},\beta_{q,\phi},\beta_{q,\theta}^{\prime},\alpha_{q,\phi}\m$ are $\pm i(c+\sqrt{ab}),\pm i(c-\sqrt{ab})$ where (see Proof of Theorem \ref{DNd}):
$$ 
\left\lbrace\begin{array}{ll}
a=d^2H_{\xi}(x_e)\cdot (\alpha_{q,\tet}^{\prime},\alpha_{q,\tet}^{\prime})\\
b=d^2H_{\xi}(x_e)\cdot (\beta_{q,\phi},\beta_{q,\phi})\\
c=d^2H_{\xi}(x_e)\cdot (\alpha_{q,\tet}^{\prime},\beta_{q,\phi}).
\end{array}\right.
$$

Second, assume that $N$ is even.
The blocks are similar to those for $N$ odd except the last block:
the block of the Hessian $d^2H_{\xi}|_{\G}(x_e)$ corresponding to the subspace $\<\alpha_{\frac{N}{2},\theta},\alpha_{\frac{N}{2},\theta}^{\prime},\alpha_{\frac{N}{2},\phi},\alpha_{\frac{N}{2},\phi}^{\prime}\m$ is a diagonal block:
$$
\left( \begin{array}{cccc} 
a^{\prime} & 0 & 0 & 0\\
0 & a & 0 & 0\\
0 & 0 & b & 0\\
0 & 0 & 0 & b^{\prime}
\end{array}\right)
\ \mbox{where}\left\lbrace\begin{array}{llll}
a^{\prime}=d^2H_{\xi}(x_e)\cdot (\alpha_{\frac{N}{2},\theta},\alpha_{\frac{N}{2},\theta})\\
a=d^2H_{\xi}(x_e)\cdot (\alpha_{\frac{N}{2},\theta}^{\prime},\alpha_{\frac{N}{2},\theta}^{\prime})\\
b=d^2H_{\xi}(x_e)\cdot (\alpha_{\frac{N}{2},\phi},\alpha_{\frac{N}{2},\phi})\\
b^{\prime}=d^2H_{\xi}(x_e)\cdot (\alpha_{\frac{N}{2},\phi}^{\prime},\alpha_{\frac{N}{2},\phi}^{\prime}).
\end{array}\right.
$$
The associated block of the linearized system $L$ is:
$$
\left( \begin{array}{cccc} 
0 & 0 & 0 & b^{\prime}\\
0 & 0 & b & 0\\
0 & -a & 0 & 0\\
-a^{\prime} & 0 & 0 & 0
\end{array}\right)
$$
Eigenvalues of this block are $\pm i\sqrt{ab}$, $\pm i\sqrt{a^{\prime}b^{\prime}}$, and bifurcations occur when $ab$ or $a^{\prime}b^{\prime}$ passes through zero.

For the two cases $N$ odd and $N$ even, some tedious computations lead to the same formulae for $a$ and $b$:
$$
\begin{array}{lll}
a & = & \epsilon N [ -(N-1)\frac{u^2}{1-u^2}+\frac{1}{1-u^2}\sum_{j=2}^N\frac{\cos \varphi_j-\cos q\varphi_j}{1-\cos \varphi_j}\\[10pt]
  & - & \sum_{j=1}^N  \frac{u^4-u^2+(1-u^2)\cos q\varphi_j+(2u^4-u^2+1-(1+u^2)\cos q\varphi_j)\cos \varphi_j-(1-u^4)\cos^2 \varphi_j}{[1+u^2-(1-u^2)\cos \varphi_j]^2} ]\\[10pt]
b & = & \epsilon N [ -\sum_{j=2}^N\frac{1-\cos q\varphi_j}{1-\cos \varphi_j}  +(1-u^2)\sum_{j=1}^N  (1-\cos q\varphi_j)\frac{1-u^2-(1+u^2)\cos \varphi_j}{[1+u^2-(1-u^2)\cos \varphi_j]^2} ]
\end{array}
$$
where $\epsilon=1$ if $q\neq N/2$ and $\epsilon=2$ if $q=N/2$.  
One checks that if $q$ is the integer part of $N/2$, then for all $\theta_0$ and all $N\geq 8$ we have $ab<0$, and therefore relative equilibria in question are linearly unstable.

A study with MAPLE gives results for $N=2,\dots,7$ showing that the blocks governing stability are those considered above.
The Lyapunov stability is lost through a pitchfork bifurcation:
$D_{2h}(2R)$ bifurcates to $C_{2v}(R_m,R_m^{\prime})$ and for $N\geq 4$, $D_{Nh}(2R)$ bifurcates to $D_{\frac{N}{2}h}(4R)$.
For $N=2$, we have explicitely derived the equations of the bifurcating branches in Section \ref{lesssym}.
Following the method of Section \ref{lesssym}, one can also derive the equations of the bifurcating branches for $N\geq 4$.
\ep

Then, we add as usual two polar vortices of vorticity $\pm 1$ and obtain the following stability results.
\begin{theo}
\label{DNhp}
\
 
\noindent $\bullet$ For $N \geq 9$, $D_{Nh}(2R,2p)$ relative equilibria are linearly unstable for all co-latitude $\theta_0$.\\
$\bullet$ $D_{2h}(2R,2p)$ relative equilibria are  linearly unstable for all $\theta_0$.\\
The following relative equilibria are Lyapunov stable in the given range  and linear unstable in the complement in $]0,\pi[$:\\
$D_{4h}(2R,2p)$  $]0,0.92]$\\
$D_{6h}(2R,2p)$  $]0,0.83]$\\
$D_{8h}(2R,2p)$  $]0,0.47]$.\\
A pitchfork bifurcation occurs at the upper bound of these intervals.\\
$\bullet$ The following relative equilibria are Lyapunov stable on the first range, linearly stable on the second range and linearly unstable on the third one:\\
$D_{3h}(2R,2p)$  $]0,0.83]$, $]0.83,0.87]$, $]0.87,\pi]$\\
$D_{5h}(2R,2p)$  $]0,0.91]$, $]0.91,0.93]$, $]0.93,\pi]$\\
$D_{7h}(2R,2p)$  $]0,0.71]$, $]0.71,0.72]$, $]0.72,\pi]$.\\
Linear stability is lost through a Hamiltonian-Hopf bifurcation.
\end{theo}

\proof 
The proof is a combination of those of Theorems \ref{DNdp} and \ref{DNh}.
\ep

\rmk 
The larger is the vorticity of the North pole, the longer is the range of stability of $D_{Nh}(2R,2p)$. 
For example, if $\lambda_n=+10$ then $D_{8h}(2R,2p)$ is stable on $]0,1.2]$.

\subsection{Other relative equilibria stability results}

In this section, we compute stability of the relative equilibria found in Section \ref{lesssym}.
Relative equilibria in question consist of $4$ or $6$ vortices.
It is not necessary to get an isotypic decomposition of the symplectic slice since its dimension is either $4$ or $8$.
A study of the eigenvalues of the linearized system in the symplectic slice with MAPLE gives the following results:\\
$\bullet$ $C_{2v}(R,R')$, $C_{2v}(R,R',2p)$ and $C_{2v}(2R,2p)$ with $\lambda_n=\pm 1$ are \emph{linearly unstable} for all co-latitudes.\\ 
$\bullet$ $C_{2v}(R_m,R_m^{\prime})$ are \emph{Lyapunov stable} provided the two (+)vortices lie in the same hemisphere and \emph{linearly unstable} if this condition is not satisfied.

\subsection{Bifurcations}

The energy diagrams (see figures 8 and 9) represent the Hamiltonian (the energy) with respect to the momentum for the relative equilibria of systems of $4$ and $6$ vortices.
Polar vortices are such that $\lambda_n=\pm 1$.

A subcritical pitchfork bifurcation occurs for $\mu=1.66$, the $D_{2d}(R,R')$ relative equilibrium bifurcates to $C_{2v}(R,R')$.

A supercritical pitchfork bifurcation occurs for $\mu=3.15$, the $D_{2h}(2R)$ relative equilibrium bifurcates to $C_{2v}(R_m,R_m^{\prime})$.\\  
Note that bifurcations are pitchfork since the two bifurcating branches are exchanged by a $\Z_2$ symmetry (see Section \ref{lesssym}).

For the following relative equilibria $D_{4h}(2R)$, $D_{6h}(2R)$, $D_{8h}(2R)$, $D_{6h}(2R,2p)$, $D_{8h}(2R,2p)$,  Lyapunov stability is lost through a pitchfork bifurcation:
$D_{Nh}$ breaks to $D_{\frac{N}{2}h}$.
And it seems likely that a second pitchfork bifurcation follows immediately:
the $D_{\frac{N}{2}h}$ symmetry breaks to $C_{\frac{N}{2}v}$ .
Following the method of Section 3.2.2, one can derive the equations of the bifurcating branches.
\begin{figure}[p]
    \begin{center}
\includegraphics[width=4.6in,height=4.6in,angle=0]{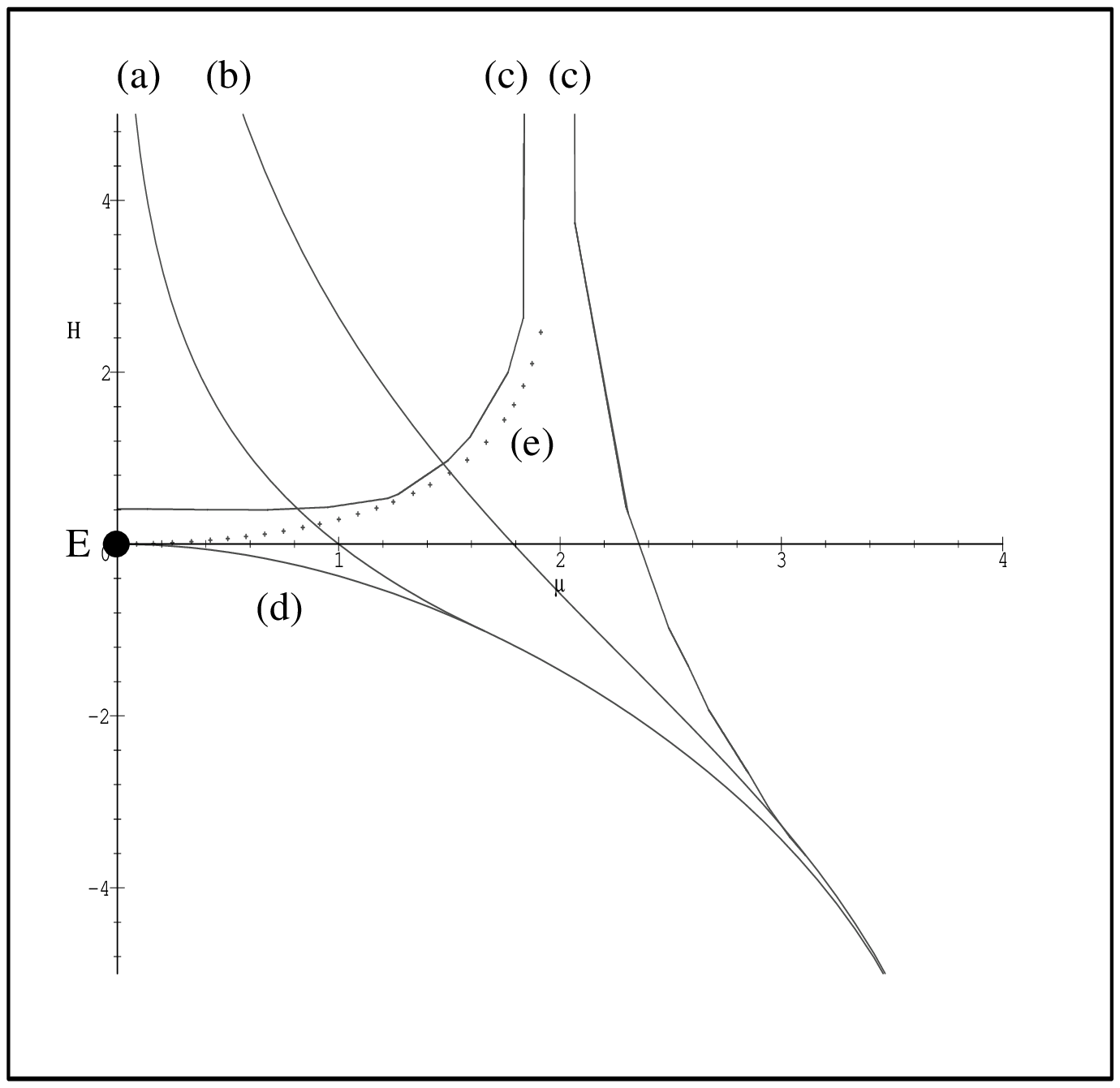}
\caption{Energy-momentum diagram for relative equilibria formed of $2$ ($+$)vortices and $2$ ($-$)vortices.
The relative equilibria are:
(a) $C_{2v}(R,R')$, (b) $D_{2h}(2R)$, (c) $C_{2v}(R_m,R_m^\prime)$, (d) $D_{2d}(R,R')$, (e) $C_{2v}(R,2p)$.
The point $E$ corresponds to the equatorial fixed equilibrium $D_{4h}(R_e)$.}
     \end{center}
\end{figure}
\begin{figure}[p]
    \begin{center}
\includegraphics[width=4.6in,height=4.6in,angle=0]{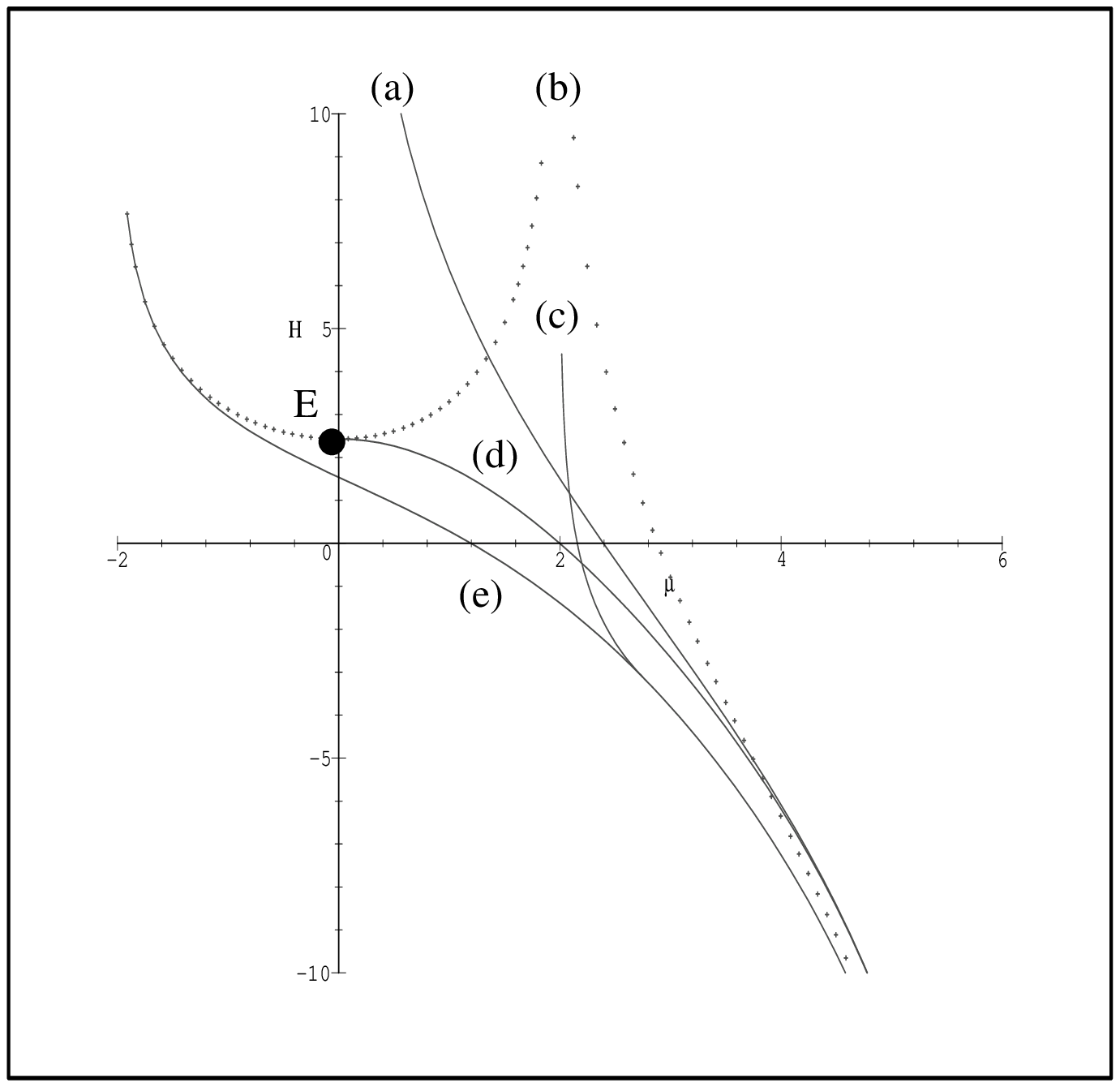}
\caption{Energy-momentum diagram for relative equilibria formed of $3$ ($+$)vortices and $3$ ($-$)vortices.
The relative equilibria are:
(a) $D_{3h}(2R)$, (b) $D_{2h}(2R,2p)$, (c) $C_{2v}(R,R',2p)$, (d) $D_{3d}(R,R')$, (e) $D_{2d}(R,R',2p)$.
The point $E$ corresponds to the equatorial fixed equilibrium $D_{6h}(R_e)$.}
     \end{center}
 \end{figure}

\subsection{Conclusion}

The larger is the number of vortices, the more unstable is the relative equilibrium.
This agrees with the following well known fact in chemistry:
the larger is the number of atoms of a molecule, the more unstable is the molecule.

Relative equilibria $D_{Nh}(2R,k_pp)$ are stable near the poles, while relative equilibria $D_{Nd}(2R,k_pp)$ are rarely stable.
One can increase the range of stability of these relative equilibria by increasing the vorticity of the polar vortices:
for $D_{Nh}(2R,k_pp)$ one needs to set the ``plus'' polar vortex and the $(+)$ring in the same hemisphere, while for $D_{Nd}(2R,k_pp)$ one needs to set the ``minus'' polar vortex with the $(+)$ring.\\
\\
\\
\textbf{Aknowledgements}\\
This work on point vortices with opposite vorticities was suggested by James Montaldi.
I would like to gratefully aknowledge and thank James Montaldi and Pascal Chossat for very helpful comments and suggestions.
A visit to the University of Warwick was supported by the MASIE network.
Finally, I thank the anonymous referees for several suggestions which have improved this paper.


\begin{thebibliography}{99}
\bibitem[Ar82] {Ar82} H. Aref, Point vortex motion with a center of symmetry. \emph{Phys. Fluids} {\bf 25} (1982), 2183-2187.
\bibitem[Ar83a] {Ar83} H. Aref, Integrable, chaotic and turbulent vortex motion in two-dimensionnal flows. \emph{Ann. Rev. Fluid Mech.} {\bf 15} (1983), 345-389.
\bibitem[Ar83b]{Ar83'} H. Aref, The equilibrium and stability of a row of point vortices. \emph{J. Fluid Mech.} {\bf 290} (1983), 167-181.
\bibitem[AV98] {Ar98} H. Aref and C. Vainchtein, Asymmetric equilibrium patterns of point vortices. \emph{Nature} {\bf 392} (1998), 769-770.
\bibitem[B77] {B77} V. Bogomolov, Dynamics of vorticity at a sphere. \emph{Fluid Dyn.} {\bf 6} (1977), 863-870.
\bibitem[CL00] {CL00} P. Chossat and R. Lauterbach, Methods in Equivariant Bifurcation and Dynamical Systems. Advanced Series in Nonlinear Dynamics {\bf 15} World Scientific (2000).
\bibitem[GE92] {GE92} L. Glasser and A. Every, Energies and spacings of point charges on a sphere. \emph{J. Phys. A} {\bf 25} (1992), 2473-2482.
\bibitem[GSS88] {GSS88} M. Golubitsky, I. Stewart, D. Schaeffer, Singularities and Groups in Bifurcation Theory, Vol.II, Applied Mathematical Sciences {\bf 69} Springer-Verlag (1988).
\bibitem[H] {H} H. Helmhotz, On integrals of the hydrodynamical equations which express vortex motion, \emph{Phil. Mag.} {\bf 33} (1867), 485-512.
\bibitem[KN98] {KN98} R. Kidambi and P. Newton, Motion of three point vortices on a sphere. \emph{Physica D} {\bf 116} (1998), 143-175.
\bibitem[KN99] {KN00} R. Kidambi and P. Newton, Collapse of three vortices on a sphere. \emph{Il Nuovo Cimento} {\bf 22} (1999), 779-791.
\bibitem[LS98] {LS98} E. Lerman  and S. Singer, Relative equilibria at singular points of the momentum map. \emph{Nonlinearity} {\bf 11} (1998), 1637-1649.
\bibitem[LR96] {LR96} D. Lewis  and T. Ratiu, Rotating $n$-gon/$kn$-gon vortex configurations. \emph{J. Nonlinear Sci.} {\bf 6} (1996), 385-414.
\bibitem[LMR01] {LMR00} C. Lim, J. Montaldi, M. Roberts, Relative equilibria of point vortices on the sphere. \emph{Physica D} {\bf 148} (2001), 97-135.
\bibitem[LMR] {LMR} C. Lim, J. Montaldi, M. Roberts, Stability of point vortices on the sphere. In preparation. 
\bibitem[Ma92] {Ma92} J. Marsden, Lectures on Mechanics. LMS Lecture Note Series {\bf 174} Cambridge University Press (1992).
\bibitem[MPS99] {MPS99} J. Marsden, S. Pekarsky and S. Shkoller, Stability of relative equilibria of point vortices on a sphere and symplectic integrators. \emph{Nuovo Cimento} {\bf 22} (1999). 
\bibitem[Mi71] {Mi71} L. Michel, Points critiques des fonctions invariantes sur une $G$-vari\'{e}t\'{e}. \emph{CR Acad. Sci. Paris} {\bf 272} (1971), 433-436.
\bibitem[Mo97] {Mo97} J. Montaldi, Persistence and stability of relative equilibria. \emph{Nonlinearity} {\bf 10} (1997), 449-466.
\bibitem[MoR99] {MoR99} J. Montaldi, M. Roberts, Relative equilibria of molecules. \emph{J. Nonlinear Sci.} {\bf 9} (1999), 53-88.
\bibitem[MoR] {MoR} J. Montaldi, M. Roberts, A note on semisymplectic actions of Lie groups. To appear in \emph{CR Acad. Sci. Paris}.
\bibitem[OR99] {OR99} J-P. Ortega, T. Ratiu, Stability of Hamiltonian relative equilibria. \emph{Nonlinearity} {\bf 12} (1999), 693-720.
\bibitem[P79] {P79} R. Palais, Principle of symmetric criticality. \emph{Comm. Math. Phys.} {\bf 69} (1979), 19-30.
\bibitem[Pa92] {Pa92} G. Patrick, Relative equilibria in Hamiltonian systems: the dynamic interpretation of nonlinear stability on a reduced phase space. \emph{J. Geom. Phys.} {\bf 9} (1992), 111-119 
\bibitem[PM98] {PM98} S. Pekarsky and J. Marsden, Point vortices on a sphere: Stability of relative equilibria. \emph{J. Math. Phys.} {\bf 39} (1998), 5894-5907.
\bibitem[PD93] {PD93} L. Polvani and D. Dritshel, Wave and vortex dynamics on the surface of a sphere. \textit{J. Fluid Mech.} {\bf 255} (1993), 35-64
\bibitem[Sa92] {Sa92} P. Saffman, Vortex Dynamics. Cambridge University Press (1992).
\bibitem[Se78] {Se78} J-P. Serre, Representations lin\'{e}aires des groupes finis. Hermann, 1978.
\bibitem[Si91] {Si91} J. Simo, D. Lewis and J. Marsden, Stability of relative equilibria. Part i: the reduced energy momentum method. \emph{Arch. Rat. Meca. Anal.} {\bf 115} (1991), 15-59.
\bibitem[To] {To} T. Tokieda, Tourbillons dansants. Preprint CRM Montreal.

\end{thebibliography}
\end{document}